\documentclass[aop,seceqn,MSNbibl,citesort,dvips]{arximspdf}
\usepackage{graphicx}

\doi{10.1214/10-AOP579}
\volume{39}
\issue{4}
\pubyear{2011}
\firstpage{1243}
\lastpage{1285}

\makeatletter

\newtheorem{theorem}{Theorem}[section]
\newtheorem{cor}[theorem]{Corollary}
\newtheorem{lem}[theorem]{Lemma}
\newtheorem{clm}[theorem]{Claim}
\newtheorem{prop}[theorem]{Proposition}

\newproclaim{rem}[theorem]{Remark}

\newcommand{\vphiD}{\vphi^{-1\prime}_D}

\newcommand{\inter}{\mathrm{INT}}
\newcommand{\ninter}{\overline{\mathrm{INT}}}
\newcommand{\exit}{\mathrm{EXIT}}

\newcommand{\N}{\mathbb{N}}
\newcommand{\Z}{\mathbb{Z}}
\newcommand{\C}{\mathbb{C}}

\newcommand{\eps}{\varepsilon}

\newcommand{\E}{\mathbb{E}}

\newcommand{\ov}[1]{\overline{#1}}

\newcommand{\p}{\partial}
\newcommand{\SLE}{\mathrm{SLE}}
\newcommand{\U}{\mathbb{U}}
\newcommand{\vphi}{\varphi}

\newcommand{\around}{\circlearrowleft}

\newcommand{\Ree}{\operatorname{Re}}
\newcommand{\Imm}{\operatorname{Im}}

\newcommand{\Ee}{\mathcal{E}}

\newcommand{\D}{\mathfrak{D}}

\newcommand{\F}{\mathcal{F}}

\newcommand{\T}{\Theta}

\newcommand{\QL}{\mathcal{QL}}
\newcommand{\X}{\mathcal{X}}
\newcommand{\diam}{\operatorname{diam}}
\newcommand{\dist}{\operatorname{dist}}
\newcommand{\Y}{\mathcal{Y}}
\newcommand{\A}{\mathcal{A}}

\makeatother

\begin{document}
\begin{frontmatter}

\title{Loop-erased random walk and Poisson kernel on~planar graphs}
\runtitle{LERW on planar graphs}

\begin{aug}
\author[A]{\fnms{Ariel} \snm{Yadin}\corref{}\ead[label=e1]{a.yadin@statslab.cam.ac.uk}}
and
\author[B]{\fnms{Amir} \snm{Yehudayoff}\ead[label=e2]{amir.yehudayoff@gmail.com}}
\runauthor{A. Yadin and A. Yehudayoff}
\affiliation{University of Cambridge and Technion---IIT}
\address[A]{Statistical Laboratory\\
University of Cambridge\\
DPMMS, 3 Wilberforce Road\\
Cambridge, CB3 0WB\\
United Kingdom\\
\printead{e1}}
\address[B]{Department of Mathematics\\
Technion---Israel Institute of Technology\\
Haifa 32000\\
Israel\\
\printead{e2}}
\end{aug}

\received{\smonth{10} \syear{2008}}

%
\begin{abstract}
Lawler, Schramm and Werner showed that the scaling limit of the
loop-erased random walk on $\Z^2$ is $\SLE_2$. We consider scaling
limits of the loop-erasure of random walks on other planar graphs
(graphs embedded into $\C$ so that edges do not cross one another).
We show that if the scaling limit of the random walk is planar
Brownian motion, then the scaling limit of its loop-erasure is
$\SLE_2$. Our main contribution is showing that for such graphs,
the discrete Poisson kernel can be approximated by the continuous
one.

One example is the infinite component of super-critical percolation on
$\Z^2$.
Berger and Biskup showed that the scaling limit of the random walk on
this graph
is planar Brownian motion.
Our results imply that the scaling limit of the
loop-erased random walk on the super-critical percolation cluster
is $\SLE_2$.
\end{abstract}

%
\begin{keyword}[class=AMS]
\kwd{60F17}
\kwd{60J99}
\kwd{60K35}.
\end{keyword}
\begin{keyword}
\kwd{Loop-erased random walk}
\kwd{Schramm--Loewner evolution}
\kwd{Poisson kernel}
\kwd{planar graphs}.
\end{keyword}

\end{frontmatter}

\section{Introduction}\label{sec1}

Let $G$ be a graph. The \textit{loop-erased random walk} or LERW on
$G$ is obtained by performing a random walk on $G$,
and then erasing the loops in the random walk path in chronological
order. The
resulting path is a self-avoiding path in the graph $G$,
starting and ending at the same points as the random walk.
LERW was invented by Lawler in~\cite{LERWdfn} as a natural measure on
self-avoiding
paths. It was studied extensively on the graphs $\Z^d$.
In dimensions $d \geq4$, the scaling limit is known to be Brownian
motion (see~\cite{LawlerSurvey}).
In dimension $d = 3$, Kozma proved that the scaling limit exists
and that the limit is invariant under rotations and dilations (see
\cite{Kozma}).
In order to study the case $d = 2$, in~\cite{SchSLE} Schramm
introduced a one-parameter family of
random continuous curves, known as Schramm--Loewner evolution or $\SLE
_\kappa$.
In~\cite{LSW} Lawler, Schramm and Werner proved that the scaling limit
of LERW on $\Z^2$ is $\SLE_2$.
Their result also holds for other two-dimensional lattices.
Many other processes in statistical mechanics have been shown to
converge to
$\SLE_\kappa$ for other values of $\kappa$.

In this paper, we focus on the scaling limit of LERW on planar graphs,
not necessarily
lattices. A planar graph is a graph embedded into the complex plane so
that edges do not intersect each other;
a precise definition is provided in Section~\ref{sec: def and not}.
We allow weighted and directed graphs, but require them to be irreducible;
that is, any two points are connected by a path in the graph.

Our main result, Theorem~\ref{thm: main thm}, is a generalization of
\cite{LSW}.
Let $G$ be an irreducible graph,
and let $f \dvtx G \to\C$ be an embedding of $G$ into the complex plane.
If $f(G)$ is planar (in the sense above),
and if the scaling limit of the random walk on $f(G)$ is planar
Brownian motion,
then the scaling limit of LERW on $f(G)$ is $\SLE_2$.

\begin{figure}

\includegraphics{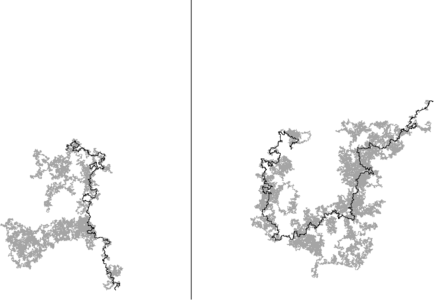}

\caption{LERW \textup{(black)} and simple random walk \textup{(gray)}
stopped on exiting the unit disc. The underlying graphs are $\Z^2$
\textup{(left)}
and the super-critical percolation cluster with parameter $0.75$
\textup{(right)}.
The mesh size is $1/600$.}\label{fig: SCP}
\end{figure}

One interesting example is the infinite component of super-critical
percolation on $\Z^2$.
That is, consider bond percolation on $\Z^2$, each bond open with
probability $p > 1/2$, all bonds independent.
Then, a.s. there exists a unique infinite connected component. In \cite
{BerBisk} Berger and Biskup
proved that a.s. the scaling limit of the random walk on this infinite
component is Brownian motion.
Together with our result, this implies that a.s. the scaling limit of
LERW on the
super-critical percolation cluster is $\SLE_2$ (see Figure~\ref{fig: SCP}).

%
%

Another example of a planar graph with random walk converging to planar
Brownian motion
is given by Lawler in~\cite{LawRWRE} (see the example following Lem\-ma~5).
For each vertex $z \in\Z^2$, define transition probabilities as follows:
the probability to go either up or down is $p(z)/2$,
and the probability to go either left or right is $(1-p(z))/2$.
Lawler proved in~\cite{LawRWRE} that if $p(z)$ are all chosen i.i.d.
such that $\mathbb{P}[p(z)=p] = \mathbb{P}[p(z)=1-p] = 1/2$,
for some $0< p < 1/2$, then a.s. the scaling limit\vadjust{\goodbreak} of the random walk
on this graph is planar Brownian motion.
Our result implies that the LERW on this graph converges to $\SLE_2$.

The main contribution of this work is Lemma~\ref{lem: H(,) are close
to lambda}, that states
that for planar graphs, the discrete Poisson kernel can be approximated
by the continuous Poisson kernel.
This result holds for any bounded domain, although the boundary
behavior can be arbitrary.
This result also holds ``pointwise,'' regardless of the local geometry
of the graph.
Perhaps it can be used to generalize other limit theorems about
processes on $\Z^2$
(such as IDLA) to more general planar graphs (e.g., the super-critical
percolation cluster).

\subsection{Definitions and notation}
\label{sec: def and not}

For any $v,u \in\C$, denote $[v,u] = \{ (1-t)v + t u \dvtx0 \leq t
\leq1 \} $. 

\subsubsection*{Planar-irreducible graphs}
Let $G = (V,E)$ be a directed weighted graph;
that is, $E\dvtx V \times V \to[0,\infty)$.
We write $(v,u) \in E$, if $E(v,u)> 0$.
Let $o \in V$ be a fixed vertex.
Let $f\dvtx V \to\C$ be an embedding of $G$ in the complex plane such that:
\begin{longlist}[(1)]
\item[(1)] $f(o) = 0$.
\item[(2)] The embedding of $G$ in $\C$ is a ``planar'' graph;
that is, for every two edges $(v,u) , (v',u') \in E$
such that $\{v,u\} \cap\{v' , u'\} =
\varnothing$,
$[f(v),f(u)] \cap[f(v'),f(u')] = \varnothing$.
\item[(3)] For every compact set $K \subset\C$, the number of vertices
$v \in V$ such that $f(v) \in K$ is finite.
\end{longlist}
We think of the graph $G$ as its embedding in $\C$.
For $\delta> 0$, let $G_\delta= (V_\delta,E_\delta)$ be the graph
defined by
\[
V_\delta= \{ \delta f(v) \dvtx v \in V \}
\quad\mbox{and}\quad E_\delta(\delta f(v), \delta f(u)) = E(v,u);
\]
that is, $G_\delta$ is the embedding of $G$ in $\C$ scaled by a
factor of $\delta$.

We assume that $\sum_{u \in V} E(v,u) < \infty$ for every $v \in V$.
Let $P\dvtx V \times V \to[0,1]$ be
\[
P(v,u) = \frac{E(v,u)}{\sum_{w \in V} E(v,w)} .
\]
%
We call the Markov chain induced on $V_\delta$ by $P$ the \textit
{natural random walk on $G_{\delta}$}.
We assume that the natural random walk is \textit{irreducible};
that is, for every $v,u \in V$, there exists $n \in\N$ such that
$P^n(v,u) > 0$.

We call a graph $G$ that satisfies all the above properties a \textit
{planar-irreducible} graph.
For the remainder of this paper we consider only planar-irreducible graphs.

\subsubsection*{Loop erasure}
Let $x(0), x(1), \ldots, x(n)$ be $n+1$ vertices in $G_{\delta}$.
Define $x[0,n]$ as the linear interpolation of $(x(0), \ldots, x(n))$;
that is, for $t \in[0,n]$, set
\[
x(t) = \bigl(1 - (t-\lfloor t \rfloor)\bigr) x(\lfloor t \rfloor) +
(t -
\lfloor t \rfloor) x(\lfloor t \rfloor+ 1) .\vadjust{\goodbreak}
\]
Define the \textit{loop-erasure} of $x(\cdot)$ as the self-avoiding
sequence induced by erasing loops in chronological order; that is,
the loop-erasure of $x(\cdot)$ is the sequence $y(\cdot)$
that is defined inductively as follows: $y(0) = x(0)$,
and $y(k+1)$ is defined using $y(k)$ as $y(k+1) = x(T+1)$, where
$T = \max\{ \ell\leq n \dvtx x(\ell) = y(k) \} $ [the loop-erasure
ends once $y(k) = x(n)$].

A \textit{path} from $v$ to $u$ in $G_\delta$ is a sequence $v=x(0),x(1),
\ldots,x(n)=u$ such that $(x(j),x(j+1)) \in E_{\delta}$ for all $j$.
The \textit{reversal} of the path
$x(\cdot)$ is the sequence $x(n),x(n-1), \ldots,x(0)$.
The reversal of a path is not necessarily a path.

%

\subsubsection*{Domains}
Denote by $\U$ the open unit disc in $\C$.
Let $D \subsetneqq\C$ be a simply connected domain such that $0 \in D$.
Define $V_\delta(D)$ as the set of vertices $z \in V_\delta\cap D$
such that there is a path from $0$ to $z$ in $G_\delta$.
Define
\[
\p V_\delta(D) = \{ (v,u) \dvtx(v,u) \in E_\delta, v \in V_\delta
(D) , [v,u] \cap\p D \neq\varnothing\} ,
\]
the ``boundary'' of $G_{\delta}$ in $D$.
Denote by $\vphi_D\dvtx D \to\U$ the unique conformal map onto
the unit disc such that $\vphi_D(0) = 0$ and $\vphi_D'(0) > 0$.
Define the \textit{inner radius of $D$} as
$\operatorname{rad}(D) = \sup\{ R \geq0 \dvtx R \cdot\U\subseteq
D \} $.

Throughout this paper, we work with a fixed domain and its sub-domains.
Fix a specific bounded domain $\mathbf{D} \subsetneqq\C$ such that
$\operatorname{rad}(\mathbf{D}) > 1/2$
(one can think of $\mathbf D$ as~$\U$). Denote
\[
\D= \{ D \subseteq{\mathbf D} \dvtx D \mbox{ simply connected
domain} , \operatorname{rad}(D) > 1/2 \} .
\]

\subsubsection*{SLE}
Radial $\SLE_{\kappa}$ in $\U$ can be described as follows
(for more details see, e.g.,~\cite{LawlerBook,LSW,RS01,SchSLE,Werner}).
Let $\gamma$ be a simple curve from $\p\U$ to $0$.
Parameterize $\gamma$ so that $g'_t(0) = e^t$,
where $g_t$ is the unique conformal map
mapping $\U\setminus\gamma[0,t]$ onto $\U$
with $g_t(0) = 0$ and $g'_t(0) > 0$.
It is known that the limit $W(t) = \lim_{z \to\gamma(t)} g_t(z)$ exists,
where $z$ tends to $\gamma(t)$ from within $\U\setminus\gamma[0,t]$.
In addition, $W\dvtx[0,\infty) \to\p\U$ is a continuous function,
and the Loewner differential equation is satisfied
\[
\p_t g_t(z) = g_t(z) \frac{W(t) + g_t(z)}{W(t) - g_t(z)}
\]
and $g_0(z) = z$.
The function $W(\cdot)$ is called the \textit{driving function} of
$\gamma$.

Taking $W(t) = e^{i B(\kappa t)}$,
where $B(\cdot)$ is a one-dimensional Brownian motion
(started uniformly on $[0,2\pi]$),
one can solve the Loewner differential equation,
obtaining a family of conformal maps $g_t$.
It turns out that for $\kappa\leq4$,
the curve $\gamma$ obtained from the driving function $W$
(defined as $\gamma(0) = W(0)$ and $\gamma(0,t] = \U\setminus
g^{-1}_t(\U)$)
is indeed a simple curve from $\p\U$ to $0$ (see~\cite{RS01}).
The curve $\gamma$ is called the $\SLE_\kappa$ path.

\subsubsection*{Weak convergence}
We define weak convergence using one of several equivalent definitions
(see Chapter III in~\cite{Shiry}, e.g.).
Let $\alpha,\beta\dvtx[0,1] \to\U$ be two continuous curves.
Let $\Phi$ be the set of continuous nondecreasing maps $\phi\dvtx[0,1]
\to[0,1]$.
We say that $\alpha$ and $\beta$ are equivalent\vadjust{\goodbreak} if $\alpha= \beta
\circ\phi$
for some $\phi\in\Phi$.
Let $\mathcal{C}$ be the set of all equivalence classes under this relation.
Define
$\varrho(\alpha,\beta) = \inf_{\phi\in\Phi} \sup_{t \in[0,1]}
\vert\alpha(t) - \beta(\phi(t))\vert$.

It is known that $\varrho(\cdot,\cdot)$ is a metric on $\mathcal{C}$.
Let $\Sigma$ be the Borel $\sigma$-algebra generated by the open sets
of $\varrho$.
Let $\mu$ be a probability measure on $(\mathcal{C},\Sigma)$.
We say that $A \in\Sigma$ is \textit{$\mu$-continuous}, if $\mu(\p
A) = 0$, where $\p A$ is the boundary of $A$.

Let $\{\mu_n\}$ be a sequence of probability measures on
$(\mathcal{C},\Sigma)$.
We say that $\{\mu_n\}$ \textit{converges weakly} to $\mu$,
if for all $\mu$-continuous events $A \in\Sigma$,
it holds that $\mu_n(A)$ converges to
$\mu(A)$.

\subsubsection*{Poisson kernel}
Let $D \in\D$.
For $a \in V_\delta(D)$ and $b \in V_\delta(D) \cup\p
V_\delta(D)$, define $H(a,b) = H^{(\delta)}(a,b; D)$ to be the
probability that a natural random walk on $G_\delta$, started at
$a$ and stopped on exiting $D$, visits $b$. That is,
\[
H(a,b) = \cases{
\mathbb{P} [\exists0 \leq k \leq\tau\dvtx S(k) = b ],
&\quad$b \in V_\delta(D)$, \cr
\mathbb{P} \bigl[\bigl(S(\tau-1) , S(\tau)\bigr) = b \bigr],
&\quad$b \in\p
V_\delta(D)$,}
\]
where $S(\cdot)$ is a natural random walk on $G_{\delta}$ started at $a$,
and $\tau$ is the exit time of $S(\cdot)$ from $D$.
We sometimes denote the segment
$(S(\tau-1) , S(\tau))$ by $S(\tau)$; for example, instead of
$(S(\tau-1),S(\tau)) = b$ we write
$S(\tau) = b$, and for a set $J \subseteq\p D$, we write
$S(\tau) \in J$ instead of writing $[S(\tau-1) , S(\tau)] \cap J
\neq\varnothing$.

Let $e = (v,u) \in\p V_\delta(D)$. Let $\tilde{e} \in\p D$
be the ``first'' point on the $[v,u]$ that is not in $D$;
that is, let $s = \inf\{ 0 \leq t \leq1 \dvtx(1-t) v + t u \notin D
\} $, and let $\tilde{e} = (1-s)v + s u$.
Define 
$\vphi(e) = \lim_{t \to s^-} \vphi((1-t)v+tu)$.

For $a \in V_\delta(D)$ and $b \in V_\delta(D) \cup\p
V_\delta(D)$, define the \textit{Poisson kernel}
\[
\lambda(a,b) = \lambda(a,b; D) = \frac{1-
\vert\vphi(a)\vert^2}{\vert\vphi(a) - \vphi(b) \vert^2} .
\]

If $B(\cdot)$ is a planar Brownian motion started at $x \in\U$,
$\tau$ is the exit time of $B(\cdot)$ from $\U$, and $J$ is a Borel
subset of $\p\U$, then
%
%
\begin{equation}
\label{eqn: harmonic measure}
\mathbb{P}_x [B(\tau) \in J ] = \int_J \lambda(x,
\zeta; \U) \,d \zeta,
\end{equation}
where $d \zeta$ is the uniform measure on $\p\U$
(see Chapter 3 of~\cite{PeresBM}).

\subsubsection*{Complex analysis}
Throughout the proofs we will make repeated use of three classical
theorems in the theory of analytic and conformal maps: the Schwarz
lemma, the Koebe distortion theorem
and the Koebe $1/4$ theorem. These can be found in \cite
{Conway} or~\cite{Pom}.

\subsection{Main results}

Let $G$ be a planar-irreducible graph.
Let $\nu_\delta$ be the law of the natural random walk on
$G_\delta$ started at $0$ and stopped on exiting $\U$. Let
$\mu_{\delta}$ be the law of the loop-erasure of the reversal of the
natural random walk on $G_{\delta}$ started at $0$ and stopped on
exiting $\U$.
\begin{theorem} \label{thm: main thm}
Let $\{\delta_n\}$ be a sequence converging to $0$.
If $\nu_{\delta_n}$ converges weakly to the
law of planar Brownian motion started at $0$ and
stopped on exiting $\U$, then $\mu_{\delta_n}$
converges weakly to the law of radial $\SLE_2$ in $\U$ started
uniformly on $\p\U$.
\end{theorem}

The proof of Theorem~\ref{thm: main thm} is given in Section \ref
{sec: conv of lerw}.
A key ingredient in the proof is the following lemma,
that shows that the discrete Poisson kernel can be approximated by the
continuous one
(its proof is given in Section~\ref{sec: con to poison ker}).

\begin{lem} \label{lem: H(,) are close to lambda}
For all $\eps, \alpha> 0$, there exists $\delta_0$ such that for all
$0 < \delta< \delta_0$ the following holds:

Let $D \in\D$, let $a \in V_\delta(D)$ be such that $\vert\vphi
_D(a)\vert\leq1-\eps$, and let
$b \in\p V_{\delta}(D)$. Then,
\[
\biggl\vert\frac{ H^{(\delta)}(a,b;D)}{ H^{(\delta)}(0,b;D) } -
\lambda
(a,b;D) \biggr\vert\leq\alpha.
\]
\end{lem}

Lemma~\ref{lem: H(,) are close to lambda} holds for
all graphs that are planar, irreducible and such that the scaling limit
of the
random walk on them is planar Brownian motion.
The question arises whether a similar result holds in ``higher
dimensions.''
The answer is negative. For $d>2$, one can construct a subgraph of $\Z
^d$ such that Lemma
\ref{lem: H(,) are close to lambda} does not hold for it. The idea is
to disconnect
one-dimensional subsets, leaving only one edge connecting them to the
rest of $\Z^d$.
This can be done in a way so that the random walk will still converge
to $d$-dimensional
Brownian motion, but for points in these sets the discrete Poisson
kernel will be far from the continuous one.

One can also ask whether Lemma~\ref{lem: H(,) are close to lambda} can
be generalized to
nonplanar graphs. The answer is again negative. Consider the
underlying graph of
the following Markov chain. Toss a coin; if it comes out heads,
run a simple random walk on $\delta\Z^2$
conditioned to exit the unit disc in the upper half plane, and if the
coin comes out tails,
run a simple random walk on $\delta\Z^2$ conditioned to exit the unit
disc in the lower half
plane. This Markov chain converges to planar Brownian motion, but the
underlying graph
is not planar. In this example, for any point other than $0$,
the discrete Poisson kernel is supported only on one half of the unit disc
(and so is far from the continuous one).

The proof of Theorem~\ref{thm: main thm} mainly follows the proof of
Lawler, Schramm and Werner in~\cite{LSW}.
To understand the new ideas in our paper, let us first give a very
brief overview of the argument in~\cite{LSW}.
Denote by $\gamma$ the loop-erasure of the reversal of the natural
random walk,
and let $W$ be the driving function of $\gamma$ given by Loewner's thoery.

The first step is to show that $W$ converges to Brownian motion on $\p
\U$.
A key ingredient in this step is showing that the discrete Poisson
kernel can be approximated by the continuous
Poisson kernel (see Lemma~\ref{lem: H(,) are close to lambda} above).
The proof of the convergence\vadjust{\goodbreak} of the Poisson kernel in~\cite{LSW} is
based on lattice properties,
whereas the proof here uses converges to planar Brownian motion from
only one vertex, namely $0$,
and the planarity of the graph.

The second step of the proof is using a compactness argument to
conclude a stronger type of convergence.
As in~\cite{LSW}, we show that the laws given by $\gamma$ are \textit{tight}
(see definition in Section~\ref{sec: tight} below).
The proof of tightness in~\cite{LSW} uses a ``natural'' family of
compact sets.
In our setting, it is not necessarily true that $\gamma$ belongs to
one of these
compact sets with high probability (and so the argument of~\cite{LSW} fails).
To overcome this difficulty, we define a ``weaker'' notion of tightness,
which we are able to use to conclude the proof.

We now discuss the first step, the proof of Lemma~\ref{lem: H(,)
are close to lambda}, in more detail. Let $a$ be a vertex in $\U$,
and let $b$ be an edge on $\p\U$ (in fact, we need to consider
arbitrary $D \in\D$, but we ignore this here). The intuition
behind Lemma~\ref{lem: H(,) are close to lambda} is that two
independent planar Brownian motions, started at $0$ and at the
vertex $a$, conditioned on exiting $\U$ at a small interval around
$b$, intersect each other with high probability. Intuitively,
this should give us a way to couple a random started at $0$ and a
random walk
started at the vertex $a$ (conditioned on exiting $\U$ at a small
interval around $b$), so that they will both exit $\U$ at the same
point with high probability.
There are several obstacles in this argument:
first, we are not able to provide such a coupling,
and we overcome this difficulty using harmonic functions.
Second, we are not given a priori any
information on the random walk starting at the vertex $a$. Third,
we also need to consider the case where the two walks do not
intersect. Finally, we are interested in what happens at a specific
edge $b$, and not in its \textit{neighborhood} (the local geometry
around $b$ can be almost arbitrary). The main properties of $G$
that allow us to overcome these obstacles are its planarity and
the weak convergence of the random walk started at $0$ to planar
Brownian motion.

\section{Preliminaries}

Let $D \in\D$.
For $z \in V_\delta(D)$, let $S_z(\cdot)$ be a natural random
walk on $G_{\delta}$ started at $z$.
Let $\tau^{(z)}_{D}$ be the exit time of $S_z(\cdot)$ from $D$.
When $D$ is clear, we omit the subscript from $\tau^{(z)}_D$ and use
$\tau^{(z)}$.
For $U \subset D$, define\looseness=-1
\[
\Theta_z(U) = \Theta^D_z(U) = \min\bigl\{ 0 \leq t \leq\tau^{(z)}
\dvtx
S_z(t) \in U \bigr\}.
\]\looseness=0
%

For a path $\gamma[T_1,T_2]$ in $D$, denote by $\vphi_D \circ\gamma
[T_1,T_2]$ the path in $\U$
that is the image of $\gamma[T_1,T_2]$ under the map $\vphi_D$.

\subsection{Encompassing a point}

For $r > 0$ and $z \in\C$, denote
$\rho(z,r) = \{ \zeta\in\C\dvtx\vert\zeta-z\vert< r \} $,
the disc of radius $r$ centered at $z$.

\subsubsection*{Crossing a rectangle}
Let $z_1,z_2 \in\C$ and $r>0$.
Define $\square(z_1,z_2,r)$ as the $4r$ by $4r+|z_2-z_1|$
open rectangle around the interval $[z_1,z_2]$; more precisely,
define $\square(z_1,z_2,r)$ as the interior of the convex hull of the
four points
$z_1-2r(u+v)$, $z_1-2r(u-v)$, $z_2+2r(u+v)$ and
$z_2+2r(u-v)$, where
$u = \frac{z_2-z_1}{\vert z_2-z_1\vert}$ and $v = u \cdot i$.\vadjust{\goodbreak}

Let $\gamma\dvtx[T_1,T_2] \to\C$ be a curve.
Let
$t_1 = \inf\{t \geq T_1 \dvtx\gamma(t) \in\rho(z_1,r)
\}$
and
$t_2 = \inf\{t \geq T_1 \dvtx\gamma(t) \in\rho(z_2,r)
\}$.
We say that
$\gamma[T_1,T_2]$ \textit{crosses} $\square(z_1,z_2,r)$, if $t_1 < t_2
\leq T_2$ and
$\gamma[t_1,t_2] \subset\square(z_1,z_2,r)$.

\subsubsection*{Encompassing a point}
Let $z \in\C$ and $r > 0$.
Define $z_1,\ldots,z_5 \in\C$ to be the following five points:
let $r' = r/20$, let
$z_1 = z - 8r' - 4r'i$, let $z_2 = z + 4r' - 4r'i$, let
$z_3 = z + 4r' + 4r'i$, let $z_4 = z - 4r' + 4r'i$ and let $z_5 = z -
4r' - 8r'i$.

We say that $\gamma[T_1,T_2]$ \textit{$r$-encompasses} $z$, denoted
$\gamma[T_1,T_2]
\around^{(r)} z$, if $\gamma[T_1,T_2]$ crosses all rectangles
$\square(z_1,z_2,r')$, $\square(z_2,z_3,r')$, $\square(z_3,z_4,r')$,
$\square(z_4,z_5,r')$.

If $\gamma[T_1,T_2] \around^{(r)} z$, then any path from $z$ to
infinity must intersect $\gamma[T_1,T_2]$;
that is, $z$ does not belong to the unique unbounded component of $\C
\setminus\gamma[T_1,T_2]$.
Also, if $\gamma[T_1,T_2] \around^{(r)} z$, there exist $\tau_1 <
\tau_2 \leq T_2$ such that
$\gamma[\tau_1,\tau_2] \around^{(r)} z$ and $\gamma[\tau_1,\tau
_2] \subset\rho(z,r)$.

\subsection{Compactness of $\D$}

Let $D \in\D$. We bound the derivative of $\vphi^{-1}_D$ at $0$.
Using the Schwarz lemma, 
since ${\vphi^{-1}_D}(0) = 0$,
we have $\operatorname{rad}(D)/ | \vphiD(0) | \leq1$.
Since $\operatorname{rad}(D) > 1/2$, we have $|\vphiD(0)| > 1/2$.
Using the Schwarz lemma again,
we have $| \vphiD(0) | \leq C'$, for $C' = \sup\{|x| \dvtx x
\in{\mathbf D} \}$.
Thus, there exists a constant $c = c(\mathbf{D}) > 0$ such that
%
%
\begin{equation} \label{eqn: der of vphi is bounded at 0}
c \leq| \vphiD(0) | \leq c^{-1}.
\end{equation}

Let $\eps> 0$.
Every map $\vphi^{-1}_D$, for $D \in\D$, can be thought of as a
continuous map on the compact domain
$K = \{ \xi\in\U\dvtx\vert\xi\vert\leq1-\eps\}$.
The set of maps $\{ \vphi^{-1}_D \}_{D \in\D}$ is pointwise
relatively compact.
Let $z \in K$, then for every $z' \in K$,
%
\[
\vert\vphi^{-1}_D(z) - \vphi^{-1}_D(z') \vert\leq|
\vphiD(\zeta) | \cdot| z - z' |
\]
for some $\zeta\in K$.
By the Koebe distortion theorem 
and (\ref{eqn: der of vphi is bounded at 0}),
there exists a constant $c_1 = c_1(\mathbf{D}) > 0$ such that
$| \vphiD(\zeta) | \leq c_1 \cdot\eps^{-3}$.
Thus, $\{ \vphi^{-1}_D \}_{D \in\D}$ is equicontinuous.
Hence, by the Arzel\'{a}--Ascoli theorem, $\{ \vphi^{-1}_D \}_{D \in
\D}$ is relatively compact
(as maps on $K$).
\begin{prop} \label{prop: compact approx}
For any $\eps,\eta> 0$, there exist $\delta_0 > 0$ and a finite
family of domains $\D_{\eps,\eta}$, such that for every $D \in\D$
there exists $\tilde{D} \in\D_{\eps,\eta}$ with the following
properties:
\begin{longlist}[(2)]
\item[(1)]\hypertarget{item:com1}$\tilde{D} \subset D$.
\item[(2)]\hypertarget{item:com2} For every $a \in D$ such that
$|\vphi_D(a)| \leq1-\eps$, we
have $|\vphi_{\tilde{D}}(a)| \leq
1-\eps/2$.
\item[(3)]\hypertarget{item:com3} For every $\xi\in\p\tilde D$, we
have $|\vphi_D(\xi)| \geq
1- \eta$.
\item[(4)]\hypertarget{item:com4} For every $\xi\in\C$ such that
$|\xi| \leq1$,
we have $|\vphi_D(\vphi^{-1}_{\tilde D}(\xi)) - \xi| \leq\eta$.
\item[(5)]\hypertarget{item:com5} For every $\xi\in\C$ such that
there exists $z$ in the closure
of $\tilde D$
with $|z-\xi| \leq\delta_0$, we have
$|\vphi_D(\xi) - \vphi_{D}(z)| \leq\eta$.
\end{longlist}
\end{prop}

We call $\tilde D$ the \textit{$(\eps,\eta)$-approximation} of $D$.\vadjust{\goodbreak}
\begin{pf*}{Proof of Proposition~\ref{prop: compact approx}}
Let $\eps_1,\eps_2 > 0$ be small enough,
and let $K = \{ \xi\in\U\dvtx| \xi| \leq1 - \eps_1 \}$.
By the relative compactness of $\{ \vphi^{-1}_D \}_{D \in\D}$ (as
maps on $K$),
there exists a finite family of domains $\D'$ such that
for every $D \in\D$ there exists $D' \in\D'$ with
%
%
\begin{equation} \label{eqn: D and tilde D close - 1}
\operatorname{dist}(\vphi^{-1}_D,\vphi^{-1}_{D'}) =
\max_{x \in K} | \vphi^{-1}_{D}(x) - \vphi^{-1}_{D'}(x) |
< \eps_2 .
\end{equation}
Set $\D_{\eps,\eta}$ to be the set of $\tilde D = \vphi
^{-1}_{D'}((1-2\eps_1)\U)$ for $D' \in\D'$.

Let $D \in\D$, let $D' \in\D'$ be the closest domain to $D$ in $\D'$
and let $\tilde D = \vphi_{D'}^{-1}((1-2\eps_1)\U)$.
By (\ref{eqn: der of vphi is bounded at 0}),
and by the Koebe distortion theorem, 
for every $z \in K$,
%
%
\begin{equation} \label{eqn: bound on der of vphi - 1}\quad
\frac{\eps_1}{C} < |\vphiD(0)| \cdot\frac{1-|z|}{8} \leq
|\vphiD(z)|
\leq|\vphiD(0)| \cdot\frac{2}{(1-|z|)^3} < \frac{C}{\eps
_1^3} ,
\end{equation}
where $C = C(\mathbf{D}) > 0$ is a constant.

We prove property \hyperlink{item:com1}{(1)}.
Using (\ref{eqn: bound on der of vphi - 1}), for every $z_1 \in\U$
such that $|z_1| = 1-\eps_1$
and $z_2 \in\U$ such that $|z_2| = 1- 2\eps_1$,
%
%
\begin{equation} \label{eqn: ball around z}
|\vphi^{-1}_D(z_1) - \vphi^{-1}_D(z_2)| = |\vphiD(\xi)|
|z_1 - z_2| \geq
\frac{\eps_1^2}{C}
\end{equation}
for some $\xi\in K$.
By (\ref{eqn: D and tilde D close - 1}), for every $z \in\tilde D$,
there exists $\zeta\in\vphi^{-1}_D((1 - 2\eps_1) \U)$ such that
$|z-\zeta| < \eps_2$.
Thus, for $\eps_2 < \frac{\eps_1^2}{C}$,
we have $\tilde D \subset\vphi^{-1}_D(K) \subset D$.

We prove property \hyperlink{item:com2}{(2)}.
Let $a \in D$ be such that $|\vphi_D(a)| \leq1-\eps$.
We first show that for $\eps_1 \leq\eps/4$,
\[
\operatorname{dist}(b,\p\tilde D) \geq c \cdot\eps^2
\]
for a constant $c = c(\mathbf{D}) > 0$,
where $b = \vphi^{-1}_{D'}(\vphi_D(a))$.
Since $2 \eps_1 < \eps$, $b \in\tilde D$.
By 
the Koebe $1/4$ theorem, 
using the Koebe distortion theorem 
and since $\vphi_{\tilde D}^{-1}(x) = \vphi_{D'}^{-1} ((1-2 \eps_1) x)$,
\begin{eqnarray*}
\operatorname{dist}(b,\p\tilde D) & \geq&
\frac{( 1 - | \vphi_{\tilde D}(b) | ) \cdot
| \vphi_{\tilde D}^{-1\prime}(\vphi_{\tilde D}(b)) | }{4}
\geq\frac{( 1 - | \vphi_{\tilde D}(b) | )^2 \cdot(1-2
\eps_1) }{ C} \\
& = &\frac{( 1 - | {\vphi_{D}(a)}/({1-2\eps_1}) |
)^2 \cdot(1-2 \eps_1) }{ C}
\geq c \cdot\eps^2.
\end{eqnarray*}

Thus, $\rho(b,\eps_2) \subset\tilde D$,
for $\eps_2 < c \cdot\eps^2$.
Thus, by (\ref{eqn: D and tilde D close - 1}), $[a,b] \subset\tilde D$,
which implies, using the Koebe distortion theorem, 
%
\begin{eqnarray*}
|\vphi_{D'}(a) - \vphi_{D}(a)| &=& |\vphi_{D'}(a) - \vphi_{D'}(b)|
= |\vphi_{D'}'(\xi)| \cdot|b-a|\\
&\leq&\frac{C}{1- | \vphi_{D'}(\xi)|} \cdot\eps_2
\leq\frac{\eps_2 \cdot C}{2 \eps_1 }
\end{eqnarray*}
for some $\xi\in\tilde D$.
Thus, for $\eps_2 < \frac{\eps_1 \cdot\eps^2 }{2 C} $,
%
%
\begin{equation} \label{eqn: a far from tilde K}
|\vphi_{\tilde D}(a)| =
\frac{ | \vphi_{D'}(a) | }{ 1 - 2 \eps_1 }
\leq\frac{ 1 - \eps+ {\eps_2 \cdot C}/({2 \eps_1 }) }{ 1 - 2
\eps_1 } < 1 - \frac{\eps}{2} .\vadjust{\goodbreak}
\end{equation}

We prove property \hyperlink{item:com3}{(3)}.
Let $\xi\in\p\tilde D$.
Let $z = \vphi^{-1}_D(\vphi_{D'}(\xi))$.
By (\ref{eqn: D and tilde D close - 1}),\break $|z-\xi| < \eps_2$.
By (\ref{eqn: ball around z}),
$\rho(z,\eps_1^2/C) \subset\vphi^{-1}_D(K)$.
Thus, for $\eps_2 \leq\eps_1^2/C$, using (\ref{eqn: bound on der of
vphi - 1}),
\[
|\vphi_{D}(\xi) - \vphi_{D}(z)| \leq|\vphi_{D}'(\zeta)| \cdot
\eps_2 \leq\frac{C \eps_2}{\eps_1} \leq\eps_1
\]
for some $\zeta\in\vphi^{-1}_D(K)$.
Since $|\vphi_{D}(z)| = |\vphi_{D'}(\xi)| = 1-2\eps_1$,
\[
|\vphi_D(\xi)| \geq|\vphi_D(z)| - |\vphi_{D}(\xi) - \vphi
_{D}(z)| \geq1 - 3\eps_1 > 1- \eta
\]
for $\eps_1 < \eta/3$.

We prove property \hyperlink{item:com4}{(4)}.
Let $\xi\in\C$ be such that $|\xi| \leq1$.
Using (\ref{eqn: D and tilde D close - 1}),
\begin{eqnarray*}
|\vphi_D(\vphi^{-1}_{\tilde D}(\xi)) - \xi| & \leq&
\bigl|\vphi_D\bigl(\vphi^{-1}_{D'}\bigl((1-2\eps_1) \xi\bigr
)\bigr) - (1-2\eps_1) \xi\bigr| +
|(1-2\eps_1) \xi- \xi| \\
& = & |\vphi_D'(\zeta)| \cdot\bigl|\vphi^{-1}_{D'}\bigl
((1-2\eps_1) \xi\bigr) -
\vphi^{-1}_{D}\bigl((1-2\eps_1) \xi\bigr)\bigr| + 2\eps_1 \\
& \leq&|\vphi_D'(\zeta)| \cdot\eps_2 + 2\eps_1
\end{eqnarray*}
for some $\zeta\in e = [\vphi^{-1}_{D'}((1-2\eps_1) \xi),\vphi
^{-1}_{D}((1-2\eps_1) \xi)]$.
Since the length of $e$ is at most $\eps_2$, and since $\eps_2 \leq
\eps_1^2/C$,
using (\ref{eqn: ball around z}), we have $e \subset\vphi^{-1}_D(K)$.
Thus, $\vphi_D(\zeta) \in K$, which implies
using (\ref{eqn: bound on der of vphi - 1}) that
$|\vphi_D'(\zeta)| \leq\frac{C}{\eps_1}$.
Choosing $\eps_2 \leq\frac{\eps_1^2}{C}$ and $3 \eps_1 \leq\eta$
the proof is complete.

We prove property \hyperlink{item:com5}{(5)}. Let $\xi\in\C$ be such that there
exists $z$ in the closure of $\tilde D$ with $|z-\xi| \leq
\delta_0$. As in property \hyperlink{item:com4}{(4)}, for $\delta_0 \leq
\eps_1^2/C$, we have $[\xi,z] \subset\vphi^{-1}_D(K)$, which
implies
\[
|\vphi_D(\xi) - \vphi_D(z)| \leq|\vphi_D'(\zeta)| \cdot\delta_0
\leq\frac{C \delta_0}{\eps_1} \leq\eta
\]
for some $\zeta\in[\xi,z]$ and $\delta_0 \leq\eta\eps_1 / C$.
\end{pf*}

\section{Preliminaries for Brownian motion}

\subsection{Brownian motion measure continuity}

\begin{prop} \label{prop: mu-contin encoppasing a point}
Let $D \subsetneqq\C$ be a simply connected domain such that $0 \in D$.
Let $\nu$ be the law of planar Brownian motion $B(\cdot)$
(started at some point in $D$ and stopped on exiting $D$).
Let $\tau$ be the exit time of $B(\cdot)$ from $D$.
Then, the following events are $\nu$-continuous:
\begin{longlist}[(2)]
\item[(1)]\hypertarget{item:encompass}
For any $r>0$ and $z \in D$ such that
$\rho(z,r) \subset D$, the event \mbox{$\{B[0,\tau] \around^{(r)} z \}$}.

\item[(2)]\hypertarget{item:hitsW}
For any disc $\rho(z,r) \subset D$, the event $\{B[0,\tau] \cap
\rho(z,r) \neq\varnothing\}$.

\item[(3)]\hypertarget{item:exitatI}
If $D = \U$, for any interval $I \subset\p\U$, the event
$ \{B(\tau) \in I \}$.
\end{longlist}
\end{prop}
\begin{pf}
We use the following claim.
\begin{clm} \label{clm: Blumenthal's}
Let $U \subset D$ be an open set, and
$\tau_{\p U} = \inf\{t \geq0 \dvtx B(t) \in\p U \}$.
Then, if $U = \rho(z,r)$ or if $U = \square(z_1,z_2,r)$ for some
$z_1,z_2 \in D$,
we have $ \mathbb{P} [\tau_1 > \tau_{\p U} ] = \mathbb
{P} [\tau_2 > \tau_{\p U} ] = 0$,
where
$\tau_1 = \inf\{t \geq\tau_{\p U} \dvtx B(t) \in U \}$
and
$\tau_2 = \inf\{ t \geq\tau_{\p U} \dvtx B(t) \notin U \cup\p U
\}$.\vadjust{\goodbreak}
\end{clm}
\begin{pf}
We prove $\mathbb{P} [\tau_1 > \tau_{\p U} ] =0$. The
proof for $\tau_2$ is similar.
Let $\F(t)$ be the $\sigma$-algebra generated by $\{B(s) \dvtx
0 \leq s \leq t \}$, and
let $\F^+(t) = \bigcap_{s > t} \F(s)$.
Since
\[
\{\tau_1 = \tau_{\p U} \} = \bigcap_{n \in\N} \biggl\{
\exists0 < \eps< \frac{1}{n} \dvtx
B(\tau_{\p U} + \eps) \in U \biggr\} \in\F^+(\tau_{\p U}) ,
\]
by Blumenthal's 0--1 law and the strong Markov property (see, e.g.,
Chapter 2 in~\cite{PeresBM}),
$\mathbb{P} [\tau_1 = \tau_{\p U} \mid\F(\tau_{\p U}) ] \in
\{0,1\}$.
Since for any small enough $\eps>0$,
$\mathbb{P} [\tau_1 \leq\tau_{\p U} + \eps] \geq
\mathbb{P} [B(\tau_{\p U} + \eps) \in U ]
\geq\frac{1}{10}$,
we have $\mathbb{P} [\tau_1 > \tau_{\p U} ] = 0$.
\end{pf}

The event $\{B[0,\tau] \around^{(r)} z \}$ is the
intersection of
four events of the form
$\{ B[0,\tau]$ crosses $\square(z_j,z_{j+1},r') \}$, for appropriate
$z_1,\ldots,z_5$, and $r'$. So it suffices to prove that for any
$\square(z_1,z_2,r) \subset D$,
the event
$\{B[0,\tau] \mbox{ crosses } \square(z_1,\break z_2,r) \}$
is $\nu$-continuous.
By definition,
\[
\{B[0,\tau] \mbox{ crosses } \square(z_1,z_2,r) \}
= \{t_1 < t_2 \} \cap\{t_2 \leq\tau\} \cap
\{B[t_1,t_2] \subset\square(z_1,z_2,r) \} ,
\]
where
$t_1 = \inf\{t \geq0 \dvtx B(t) \in\rho(z_1,r) \}$ and
$t_2 = \inf\{t \geq0 \dvtx B(t) \in\rho(z_2,r) \}$.

Let
$\tau_1 = \inf\{ t \geq0 \dvtx B(t) \in\p\rho(z_1,r) \} $.
The boundary of the event $\{t_1 < t_2 \}$ is contained in
the event $\{t_1 > \tau_1 \}$.
Thus, by Claim~\ref{clm: Blumenthal's}, the boundary of $\{t_1 <
t_2 \}$ has zero $\nu$-measure.

Let $\tau_2 = \inf\{ t \geq0 \dvtx B(t) \in\p\rho(z_2, r) \} $.
The boundary of the event $\{t_2 \leq\tau\}$ is
contained in the event $\{t_2 > \tau_2 \}$.
Thus, by Claim~\ref{clm: Blumenthal's}, the boundary of $\{t_2
\leq\tau\}$ has zero $\nu$-measure.

Let
$\tau_3 = \inf\{ t_1 \leq t \leq t_2 \dvtx B(t) \in\p\square
(z_1,z_2,r) \} $ and
$\tau_4 = \inf\{ t \geq\tau_3 \dvtx B(t) \notin\square(z_1,z_2,r)
\cup\p\square(z_1,z_2,r) \}$.
The boundary of the event $\{B[t_1,t_2] \subset\square
(z_1,z_2,r) \}$ is contained
in the event $\{\tau_4 > \tau_3 \}$. Thus, by Claim \ref
{clm: Blumenthal's},
the boundary of $\{B[t_1,t_2] \subset\square(z_1,z_2,r)
\}$ has zero $\nu$-measure.

This proves property \hyperlink{item:encompass}{(1)}. A similar (simpler)
argument proves
property~\hyperlink{item:hitsW}{(2)}.
To prove property \hyperlink{item:exitatI}{(2)},
note that the measure $\nu$ is supported on curves that intersect $\p
\U$ at most at one point.
Hence, up to zero $\nu$-measure, the boundary of the event $\{
B(\tau) \in I \}$ is the event
$\{B(\tau) \in\{w,w'\} \}$, where $w$ and
$w'$ are the endpoints of $I$ in $\p\U$. Since
$\{B(\tau) \in\{w,w'\} \}$ has zero $\nu
$-measure, we are done.
\end{pf}

\subsection{Probability estimates}
\label{sec: bM}

This section contains some lemmas regarding planar Brownian
motion. Some of these lemmas may be considered ``folklore.''
For the sake of brevity, we omit the proofs.

\subsubsection*{Notation}
In the following $B(\cdot)$ is a planar Brownian motion.
For $x \in\U$, $\mathbb{P}_x$ is the measure of $B(\cdot)$
conditioned on $B(0)=x$.
For $r>0$, define $A(r)$ to be the annulus of inner radius $r$ and outer
radius $5r$ centered at $1$, intersected with the unit disc;
that is,
$A(r) = \{ 1+ z \dvtx r < \vert z\vert< 5r \} \cap\U$.
Also, define
$\xi(r) = 1 - 3r \in A(r)$.
Note that $\rho(\xi(r), r) \subset A(r)$
for $r < 1/25$.

The following proposition is a corollary of Theorem 3.15 in~\cite{PeresBM}.
\begin{prop}
\label{prop: prob b hits bal of gamma} Let $0 \neq x \in\U$ and let
$0 < c < \vert x\vert$. Let $\tau$ be the
exit time of $B(\cdot)$ from $\U$. Then,
\[
\mathbb{P}_x \bigl[\exists t \in[0,\tau] \dvtx\vert B(t)\vert\leq c
\bigr]
\geq\frac{1- \vert x\vert}{ -\log c } .
\]
\end{prop}
\begin{prop}
\label{prop: crossing and encompassing}
There exists $c>0$ such that the following holds:

Let $r>0$ and let $z \in\C$.
Let $T$ be the exit time
of $B(\cdot)$ from $\rho(z,r)$. Then for every $x \in\rho(z,r/2)$,
$\mathbb{P}_x [B[0,T] \around^{(r)} z ] \geq c$.
\end{prop}
\begin{prop} \label{prop: encompassing a point}
For any $0 < \eps< 1$, there exists $c > 0$ such that the
following holds:

Let $a \in\U$ be such that $\vert a\vert\leq1-\eps$.
Let $\tau$ be the exit time of $B(\cdot)$ from $\U$. Then,
$\mathbb{P}_0 [B[0,\tau] \around^{(\eps)} a ] \geq c$.
\end{prop}
\begin{lem} \label{lem: prob. of intersection in arc}
There exists $c>0$ such that the following holds:

Let $0 < r < \frac{1}{25}$, let $A = A(r)$ and $\xi= \xi(r)$.
Let $x \in A$ be such that $2r \leq\vert x-1\vert\leq4r$.
Let $T$ be the exit time of $B(\cdot)$ from $A$.
Then,
\[
\mathbb{P}_x \bigl[B[T_{\xi},T_\rho] \around^{(r)} \xi,
T_\rho< T \bigr] \geq c \cdot\frac{1-\vert x\vert}{r} \geq
\frac{c}{2} \cdot\frac{1-\vert x\vert^2}{r} ,
\]
where
$T_{\xi} = \inf\{ t> 0 \dvtx B(t) \in\rho(\xi,r/20) \} $
and
$T_\rho= \inf\{ t \geq T_{\xi} \dvtx B(t) \notin\rho(\xi,r) \} $.
\end{lem}
\begin{lem} \label{lem: CONDITIONAL prob. of intersection in arc}
There exists $c > 0$ such that the following holds:

Let $0 < \beta< \frac{1}{25\pi}$, and let $I = \{ e^{i t} \dvtx-
\pi\beta\leq t \leq\pi\beta\} $ be the interval on the unit
circle centered at $1$ of measure $\beta$. Let $\pi\beta\leq
r < \frac{1}{25}$, let $A = A(r)$ and $\xi= \xi(r)$.
Let $x \in A$ be such that $2r \leq\vert x-1\vert\leq4r$.
Let $\tau$ be the exit time of $B(\cdot)$ from $\U$, and let
$T$ be the exit time of $B(\cdot)$ from $A$. Then,
\[
\mathbb{P}_x \bigl[B[T_{\xi},T_\rho] \around^{(r)} \xi, T_\rho< T
\mid B(\tau) \in I \bigr] \geq c ,
\]
where
$T_{\xi} = \inf\{ t> 0 \dvtx B(t) \in\rho(\xi,r/20) \} $
and
$T_\rho= \inf\{ t \geq T_{\xi} \dvtx B(t) \notin\rho(\xi,r) \} $.
\end{lem}
\begin{lem} \label{lem: CONDITIONAL prob. of intersection a ball}
For every $\eta> 0$, there exists $c > 0$ such that the following holds:

Let $\beta,I,r,A,\xi,x,\tau$ and $T$ be as in Lemma~\ref{lem:
CONDITIONAL prob. of intersection in arc}.
Then,
\[
\mathbb{P}_x [T_{\xi,\eta} < T \mid B(\tau) \in I ] \geq c ,
\]
where
$T_{\xi,\eta} = \inf\{ t> 0 \dvtx B(t) \in\rho(\xi,\eta r) \} $.
\end{lem}
\begin{lem} \label{lem: BM from 0 several}
There exist $K ,c > 0$ such that the following holds:

Let $0 < \pi\beta< r < \frac{1}{2K}$, and
let $I = \{ e^{i t} \dvtx- \pi\beta\leq t \leq\pi\beta\} $ be
the interval on the unit circle centered at $1$ of measure
$\beta$. Let $\xi= \xi(r)$.
Let $\tau$ be the exit time of $B(\cdot)$ from $\U$.
Then,
\[
\mathbb{P}_0 \bigl[B[T_{\xi},\tau] \around^{(r)} \xi, \tau< T_{Kr}
\mid B(\tau) \in I \bigr] \geq c ,
\]
where
$T_{\xi} = \inf\{ t > 0 \dvtx B(t) \in\rho(\xi,r/20) \} $
and
$T_{Kr} = \inf\{ t > T_{\xi} \dvtx|B(t)-1| \geq Kr \} $.\vadjust{\goodbreak}
\end{lem}
\begin{lem} \label{lem: BM from 0 hits I+ I-}
There exist $K ,c > 0$ such that the following holds:

Let $\beta, r ,I, \xi,\tau, T_{\xi}$ and $T_{Kr}$ be as in
Lemma~\ref{lem: BM from 0 several}.
Then,
%
%
\begin{equation}
\label{eqn: lem BM I+}
\mathbb{P}_0 [T_{\xi} < \tau< T_{Kr} , B(\tau) \in I_{+} \mid
B(\tau) \in I ] \geq c
\end{equation}
and
%
%
\begin{equation}
\label{eqn: lem BM I-}
\mathbb{P}_0 [T_{\xi} < \tau< T_{Kr} , B(\tau) \in I_{-} \mid
B(\tau) \in I ] \geq c ,
\end{equation}
where $I_{+} = \{ e^{it} \dvtx\pi\beta/ 2 \leq t \leq\pi\beta\}$
and
$I_{-} = \{ e^{it} \dvtx-\pi\beta\leq t \leq- \pi\beta/2 \}$.
\end{lem}

\section{Planarity and global behavior}

\subsection{Continuity for a fixed domain}

\begin{prop} \label{prop: close points intersect - fixed a}
For all $\alpha> 0$, there exists $\eta> 0$ such that
for all $\eps> 0$, for all simply connected domains $D \subsetneqq\C
$ such that $0 \in D$,
and for all $\tilde{a} \in(1-\eps)\U$,
there exists $\delta_0 > 0$ such that for all $0 < \delta< \delta_0$
the following holds:

Let $y \in V_\delta(D) \cap\vphi^{-1}_D(\rho(\tilde a, \eta\eps))$.
Then, for every continuous curve $g$ starting in $\rho(\tilde a,\eta
\eps)$
and ending outside of $\rho(\tilde a,\eps)$,
the probability that $\vphi_D \circ S_y$ does not cross $g$
before exiting $\rho(\tilde a,\eps)$ is at most $\alpha$.
\end{prop}
\begin{pf}
Denote $\vphi= \vphi_D$.
For $x \in D$ and $r>0$,
define
\[
\tau^{(x)}(r) = \Theta_x (\vphi^{-1}(\rho(\tilde a,r) ) ) ,
\]
the time $\vphi\circ S_x$ hits $\rho(\tilde a,r)$, and define
\[
T^{(x)}(r) = \min\bigl\{ \tau^{(x)}(r/20) \leq t \leq\tau^{(x)}
\dvtx
\vphi
(S_x(t)) \notin\rho({\tilde a}, r) \bigr\}.
\]

We use the following claim and its corollary below.
\begin{clm} \label{clm: x encompasses twice}
There exists a universal constant $c > 0$ such that for all
$0 < r < \eps/40$, there exists $\delta_0 > 0$ such that for all $0 <
\delta< \delta_0$
the following holds:

There exists
$x \in V_{\delta}(D)$ such that $\vphi(x) \in\rho({\tilde a},
r/20)$ and
\[
\mathbb{P} \bigl[\vphi\circ S_x\bigl[0,T^{(x)}(r)\bigr] \around
^{(r)} {\tilde
a} , \vphi\circ S_x\bigl[T^{(x)}(r),T^{(x)}(20r)\bigr] \around^{(20r)}
{\tilde a} \bigr] \geq c .
\]
\end{clm}
\begin{pf}
Consider the event
\[
F = \bigl\{\vphi\circ S_0\bigl[\tau^{(0)}(r/20), T^{(0)}(r)\bigr]
\around
^{(r)} {\tilde a} , \vphi\circ S_0\bigl[T^{(0)}(r), T^{(0)}(20r)\bigr]
\around^{(20r)} {\tilde a} \bigr\}.
\]

Let $B(\cdot)$ be a planar Brownian motion, and let $\tau^{(B)}$ be
the exit time of $B(\cdot)$ from $\U$.
Let $\tau^{(B)}(r/20) = \inf\{0 \leq t \leq\tau^{(B)} \dvtx
B(t) \in\rho({\tilde a},r/20) \}$,
and let
$T^{(B)}(r) = \inf\{ \tau^{(B)}(r/20) \leq t \leq\tau^{(B)} \dvtx
B(t) \notin\rho({\tilde a},r) \}$
[$T^{(B)}(20r)$ is defined similarly].

By weak convergence and Proposition~\ref{prop: mu-contin encoppasing a
point}, by
the conformal invariance of Brownian motion, by the strong Markov property
and by Proposition~\ref{prop: crossing and encompassing}, for small
enough $\delta_0$,
%
%
\begin{eqnarray} \label{eqn: BM encompass twice}
\mathbb{P} [F ] & \geq&\frac{1}{2} \mathbb{P}_0
\bigl[B\bigl[\tau^{(B)}(r/20),T^{(B)}(r)\bigr] \around^{(r)}
{\tilde a}
,\nonumber\\
&&\hspace*{24pt}B\bigl[T^{(B)}(r),T^{(B)}(20r)\bigr] \around
^{(20r)} {\tilde a} \bigr]
\nonumber\\
& \geq&\frac{1}{2}
\mathbb{P}_0 \bigl[\tau^{(B)}(r/20) < \tau^{(B)} \bigr]
\nonumber\\[-8pt]\\[-8pt]
&&{} \times\inf_{\xi\in\rho({\tilde a},r/20)} \mathbb{P}_{\xi
} \bigl[B\bigl[0,T^{(B)}(r)\bigr] \around^{(r)} {\tilde a}
,\nonumber\\
&&\hspace*{75.7pt}B\bigl[T^{(B)}(r),T^{(B)}(20r)\bigr] \around
^{(20r)} {\tilde a} \bigr]
\nonumber\\
& \geq& c_1 \cdot\mathbb{P}_0 \bigl[\tau^{(B)}(r/20) < \tau^{(B)}
\bigr] ,\nonumber
\end{eqnarray}
where $c_1>0$ is a universal constant.
In addition, by the strong Markov property,
\begin{eqnarray*} 
\mathbb{P} [F ] & \leq&\mathbb{P} \bigl[\tau
^{(0)}(r/20) < \tau^{(0)} \bigr] \\
&&{} \times\max_{x}
\mathbb{P} \bigl[\vphi\circ S_x\bigl[0,T^{(x)}(r)\bigr] \around
^{(r)} {\tilde
a} , \vphi\circ S_x\bigl[T^{(x)}(r),T^{(x)}(20r)\bigr] \around^{(20r)}
{\tilde a} \bigr] ,
\end{eqnarray*}
where the supremum is over $x \in V_{\delta}(D) \cap
\vphi^{-1}(\rho({\tilde a},r/20))$. Hence,
since for small enough $\delta_0$,
\[
\mathbb{P} \bigl[\tau^{(0)}(r/20) < \tau^{(0)} \bigr] \leq2
\mathbb{P}_0 \bigl[\tau^{(B)}(r/20) < \tau^{(B)} \bigr] ,
\]
using (\ref{eqn: BM encompass twice}), there exists
$x \in V_{\delta}(D) \cap\vphi^{-1}( \rho({\tilde a}, r/20) )$ such that
\[
\mathbb{P} \bigl[\vphi\circ S_x\bigl[0,T^{(x)}(r)\bigr] \around
^{(r)} {\tilde
a} , \vphi\circ S_x\bigl[T^{(x)}(r),T^{(x)}(20r)\bigr] \around^{(20r)}
{\tilde a} \bigr] \geq c .
\]
\upqed
\end{pf}
\begin{cor} \label{cor: w encompasses}
There exists a universal constant $c > 0$ such that for all
$0 < r < \eps/40$, there exists $\delta_0 > 0$ such that for all $0 <
\delta< \delta_0$
the following holds:

For every
$w \in V_{\delta}(D)$ such that $\vphi(w) \in\rho({\tilde a}, r/20)$,
\[
\mathbb{P} \bigl[\vphi\circ S_w\bigl[0,T^{(w)}(20r)\bigr] \around^{(20r)}
{\tilde a} \bigr] \geq c .
\]
\end{cor}
%
%
\begin{pf}
We claim that there exists a set of vertices $U$ in
$V_{\delta}(D)$ 
such that every path in $G_\delta$ that starts from $w$
and reaches outside of $\vphi^{-1}(\rho(\tilde a, r))$, intersects $U$,
and such that
%
%
\begin{equation} \label{eqn: the set U}
\mathbb{P} \bigl[\vphi\circ S_u\bigl[0,T^{(u)}(20r)\bigr] \around^{(20r)}
{\tilde a} \bigr] \geq c
\end{equation}
for every $u \in U$, where $c > 0$ is the universal constant from
Claim~\ref{clm: x encompasses twice}.
This implies the corollary, since
$\mathbb{P} [\vphi\circ S_w[0,T^{(w)}(20r)] \around^{(20r)}
{\tilde a} ]$ is a convex
sum of $ \mathbb{P} [ \vphi\circ S_u[0,T^{(u)}(20r)] \around
^{(20r)} {\tilde a} ]$
for $u \in U$
(because $G$ is irreducible).

Indeed, let $U$ be the set of all vertices in $V_\delta(D) \cap\vphi
^{-1}(\rho(\tilde{a}, r))$
such that~(\ref{eqn: the set U}) holds.
Assume toward a contradiction that
there is a path $Y$ in $G_\delta$ starting from $w$ and reaching the
outside of $\vphi^{-1}(\rho(\tilde a, r))$,
such that $Y \cap U = \varnothing$. Then,
every path in $G_\delta$ whose image under $\vphi$ $r$-encompasses
$\tilde a$, must intersect $Y$.
Let $x$ be the vertex guarantied by Claim~\ref{clm: x encompasses twice}.
Then,
\begin{eqnarray*}
&&
\mathbb{P} \bigl[ \vphi\circ S_x\bigl[0,T^{(x)}(r)\bigr] \around
^{(r)} {\tilde a}
,
\vphi\circ S_x\bigl[T^{(x)}(r),T^{(x)}(20r)\bigr] \around^{(20r)}
{\tilde a} \bigr]
\\
&&\qquad\leq\sum_{y \in Y} p(y) \cdot
\mathbb{P} \bigl[\vphi\circ S_y\bigl[0,T^{(y)}(20r)\bigr] \around^{(20r)}
{\tilde a} \bigr] < c,
\end{eqnarray*}
which is a contradiction to Claim~\ref{clm: x encompasses twice}
[where $\{p(y)\}_{y \in Y}$ is a distribution on the set~$Y$].
%
\end{pf}

We continue with the proof of Proposition~\ref{prop: close points
intersect - fixed a}.
Let $c>0$ be the constant from Corollary~\ref{cor: w encompasses}.
Let $M \in\N$ be large enough so that $(1-c)^M < \alpha$.
Let $\eta> 0$ be small enough so that $500^{M+1} \eta< 1/40$.
For $j=1,2,\ldots,M$, define $r_j = 500^j \eta\eps$, and define
$F_j = \{\vphi\circ S_y[T^{(y)}(r_j), T^{(y)}(400r_j)] \around
^{(400r_j)} {\tilde a} \}$.
By the strong Markov property and by Corollary~\ref{cor: w encompasses},
since $\vphi(y) \in\rho(\tilde a, \eta\eps)$, we have
$ \mathbb{P} [F_j \mid\ov{F}_1 , \ldots, \ov{F}_{j-1} ] \geq c$
for every $j$,
which implies
\[
\mathbb{P} [\ov{F}_1 , \ldots, \ov{F}_{M} ] \leq
(1-c)^M < \alpha
\]
(here and below $\ov{E}$ is the complement of the event $E$).
Since $G$ is planar-irreducible,
the proposition follows.
\end{pf}

\subsection{Starting near the boundary}

In this section we prove the version of Lemma 5.4 in~\cite{LSW} that
is relevant to us.
Part of the proof is similar to that of~\cite{LSW},
but the setting here is more general and requires more details.
%
%
\begin{lem} \label{lem: starting near boundary}
For any $\eps,\alpha> 0$, there exist $\eta,\delta_0 > 0$ such that
for every $0 < \delta< \delta_0$ the following holds:

Let $D \in\D$, and let $x \in V_\delta(D)$ be such that $|\vphi
_D(x)| \geq1- \eta$.
Then, the probability that $S_x$ hits the set $\{y \in D \dvtx|\vphi
_D(y) - \vphi_D(x)| > \eps\}$
before exiting $D$ is at most~$\alpha$.
\end{lem}

We first prove the following proposition.
\begin{prop} \label{prop: starting near boundary}
There exists $0 < \alpha< 1$ such that for any $\eps> 0$,
there exist $\eta,\delta_0 > 0$ such that for every $0 < \delta<
\delta_0$ the following holds:

Let $D \in\D$, and let $x \in V_\delta(D)$ be such that $1 - 2\eta
\leq|\vphi_D(x)| \leq1 - \eta$.
Then, the probability that $S_x$ hits the set $\{y \in D \dvtx|\vphi
_D(y) - \vphi_D(x)| > \eps\}$
before exiting $D$ is at most $\alpha$.
\end{prop}
\begin{pf}
Let $\eta> 0$ be small enough.
By (\ref{eqn: der of vphi is bounded at 0}), by the Koebe distortion
theorem, 
and by the Koebe $1/4$ theorem, 
$\operatorname{dist}(x,\p D) \geq r_0$, where $r_0 = c \cdot\eta^2$ for
some constant $c = c(\mathbf{D}) > 0$.
Let $z \in\p D$ be a point such that $r = |z-x| = \operatorname
{dist}(x,\p D)$.
Let $x' \in{\mathbf D}$ be such that $|x'-x| < r_0/C$,
and let $z' \in{\mathbf D}$ be\vadjust{\goodbreak} such that $|z'-z| < r_0/C$,
for a large enough constant $C > 0$.
We need to consider only finitely many points $x'$ and $z'$.

Let $r' = |x'-z'|$, and let $R > 0$ be large enough so that ${\mathbf D}
\cup\rho(x',10r') \subset\frac{R}{2} \U$.
Denote $A_1 = \{\xi\in\C\dvtx|\xi-z'| \leq r'/10 \}
\cup[x',z'] \setminus\{x'\}$.
Let $Q$ be the connected component in $\C$ of $(\p\rho(z',r')) \cap
D$ that contains $x'$.
Let $A_2$ and $A_3$ be the two connected components in $\C$ of $Q
\setminus\{x'\}$.
For large enough $C$, the distance from $x'$ to $\p D$ is at least $3r'/4$.
Thus, both $A_2$ and $A_3$ are arcs of length at least $3r'/4$.
If $C$ is large enough, $D \setminus(A_1 \cup A_2 \cup A_3)$ has three
connected components in $\C$.
For $j = 1,2,3$, let $K_j$ be the connected component in $\C$ of
$D \setminus(A_1 \cup A_2 \cup A_3)$
such that $A_j \cap\p K_j = \varnothing$.
Let $\Ee_j$ be the collection of curves $\gamma\subset R \U$ such that
$\gamma$ stays in $K_j$ from the first time it first hits $\p\rho
(x',r'/2)$ until it exits $D$.
By the conformal invariance of Brownian motion,
there exists a universal constant $c_1 > 0$ such that for every $j = 1,2,3$,
we have $\mathbb{P}_{x'}[B(\cdot) \in\Ee_j] > c_1$, where $B(\cdot
)$ is a planar Brownian motion
started at $x'$.

Let $A = \{y \in D \dvtx|\vphi_D(y) - \vphi_D(x)| > \eps\}$. We show
that there exists $j' \in\{1,2,3\}$ such that $A \cap K_{j'} =
\varnothing$. Assume toward a contradiction that $A \cap K_j
\neq\varnothing$ for all~$j$. We prove for the case that $A$ intersects
both $A_1$ and $A_2$ (the proof for the other cases is similar). $A$ is
a connected set that intersects both $A_1$ and $A_2$, so we can choose
$A'$ to be a minimal connected subset of $A$ that intersects both $A_1$
and $A_2$ (minimal with respect to inclusion). Thus, either $A'$ is in
the closure of $K_3$ or $A'$ is in the closure of $K_1 \cup K_2$. We
prove for the case that $A'$ is in the closure of $K_3$ (the proof for
the other case is similar).

We show that $A \cap\rho(x',r'/2) = \varnothing$. By choosing $\eta
> 0$
to be small enough, and by the conformal invariance of Brownian motion,
the probability that a Brownian motion started at $x$ hits $A$ before
exiting $D$ can be made arbitrarily small. If $A \cap\rho(x,3r/5)
\neq\varnothing$, because $\operatorname{dist}(x,\p D) = r$ and because
$A$ is connected, the probability that a Brownian motion started at $x$
hits $A$ before exiting $D$ is at least a universal constant $c_2>0$.
This is a contradiction for a small enough $\eta$, which implies $A
\cap\rho(x,3r/5) = \varnothing$. Since $r' \leq r(1+ 2/C)$ and since
$|x-x'| \leq r_0/C$, for large enough $C$ we have that $\rho(x',r'/2)
\subset\rho(x,3r/5)$.

For a vertex $y \in V_\delta(\rho(x',r'/2))$,
define $h(y)$ as the probability that $S_y[0,\tau^{(y)}_D]$ is in $\Ee_{j'}$.
The map $h(\cdot)$ is harmonic in $V_\delta(\rho(x',r'/2))$
with respect to the law of the natural random walk on $G_\delta$.
\begin{clm}
\label{clm: there is z in rhoxr2}
There exist a universal constant $c_3 > 0$ and $\delta_0 > 0$ such
that for all
$0 < \delta< \delta_0$, there exists $y \in V_\delta(\rho(x',r_0/C))$
with $h(y) \geq c_3$.
\end{clm}
\begin{pf}
We prove for the case $j'=3$. The proof of the other cases is similar.
The event $\Ee_{3}$ contains an event $\Ee$ that is independent of
$D$; for example,
there exist $x'=z_1,z_2,\ldots,z_m \in\C$ for $m \leq10^3$ such
that $|z_{i+1}-z_i| = r'/2$, and
\[
\Ee= \{\gamma\subset R \U\dvtx\gamma\mbox{ crosses }
\square(z_i,z_{i+1},r'/100) \mbox{ for all } i \} \subset
\Ee_3 .\vadjust{\goodbreak}
\]

Let $B(\cdot)$ be a Brownian motion, and let $\tau$ be the exit time
of $B(\cdot)$ from $R \U$.
Since $\square(z_1,z_2,r'/100),\ldots,\square(z_{m-1},z_{m},r'/100)$
are $m-1$
rectangles of fixed proportions, we have
$\inf_{w \in\rho(x',r'/100) } \mathbb{P}_w [ B[0,\tau] \in\Ee] >
c_4$ for some universal constant $c_4 > 0$.
Let $T$ be the time $B(\cdot)$ hits $\rho= \rho(x',r_0/C)$.
On one hand,
\[
\mathbb{P}_0 \bigl[B[T,\tau] \in\Ee\bigr] \geq\mathbb{P}_0 [ T
< \tau] \cdot c_4 .
\]
On the other hand, using weak convergence and Proposition~\ref{prop:
mu-contin encoppasing a point},
if $\delta_0$ is small enough,
\begin{eqnarray*}
\mathbb{P}_0 \bigl[B[T,\tau] \in\Ee\bigr] & \leq&2 \mathbb{P}
\bigl[S_0\bigl[\Theta_0(\rho),\tau^{(0)}_{R \U}\bigr] \in\Ee
\bigr] \\
& \leq&4 \mathbb{P}_0 [ T < \tau] \cdot\max_{y \in V_\delta(\rho
)} \mathbb{P} \bigl[S_y\bigl[0,\tau^{(y)}_{R \U}\bigr] \in\Ee
\bigr].
\end{eqnarray*}
\upqed
\end{pf}

Let $c_3 > 0$ and let $y \in V_\delta(\rho(x',r_0/C))$ be given by
Claim~\ref{clm: there is z in rhoxr2}.
Since $h(\cdot)$ is harmonic,
there exists a path $\gamma$ from $y$ to $\p\rho(x',r'/2)$
such that $h(w) \geq h(y)$ for every $w \in\gamma$.
Since $h(\cdot)$ is nonnegative, harmonic and bounded,
\[
h(x) \geq\mathbb{P}\bigl[ S_x\bigl[0,\tau^{(x)}_{\rho
(x',r'/2)}\bigr] \cap\gamma
\neq\varnothing\bigr]
\cdot h(y).
\]
By Proposition~\ref{prop: close points intersect - fixed a},
and by choosing\vspace*{1pt} large enough $C$, we have
$\mathbb{P}[ S_x[0,\tau^{(x)}_{\rho(x',r'/2)}] \cap\gamma\neq
\varnothing] \geq1/2$.
Since every curve in $\Ee_{j'}$ does not intersect $A$,
the probability that $S_x$ hits the set $A$ before exiting $D$ is at most
$1-c_3/2 < 1$.
\end{pf}

Planarity and Proposition~\ref{prop: starting near boundary} imply a
stronger statement.
\begin{cor} \label{cor: starting near boundary}
There exists $0 < \alpha< 1$ such that for any $\eps> 0$,
there exist $\eta,\delta_0 > 0$ such that
for every $0 < \delta< \delta_0$ the following holds:

Let $D \in\D$, and let $x \in V_\delta(D)$ be such that $|\vphi
_D(x)| \geq1 - \eta$.
Then, the probability that $S_x$ hits the set $\{y \in D \dvtx|\vphi
_D(y) - \vphi_D(x)| > \eps\}$
before exiting $D$ is at most~$\alpha$.
\end{cor}
\begin{pf}
Let $\alpha,\eta,\delta_0$ be given by Proposition~\ref{prop:
starting near boundary} with $\eps/10$,
and let $0 < \delta< \delta_0$.
For $y \in V_\delta(D)$, define $f(y)$ as
the probability that $S_y$ hits $A=\{y \in D \dvtx|\vphi_D(y) - \vphi
_D(x)| > \eps\}$ before exiting $D$.
Assume toward a contradiction that $f(x) > \alpha$.
The map $f(\cdot)$ is harmonic in $V_\delta(D \setminus A)$
with respect to the law of the natural random walk on $G_\delta$.
Let $A'$ be the set of $\xi\in D$ such that $1-2\eta\leq|\vphi
_D(\xi)| \leq1- \eta$
and $|\vphi_D(\xi) - \vphi_D(x)| \leq\eps/2$.
By Proposition~\ref{prop: starting near boundary}, $f(y) \leq\alpha$
for all $y \in V_\delta(A')$.
Thus, there exists a path $\gamma$ from $x$ to the set $A$ in
$V_\delta(D)$ that does not intersect
$A'$ such that $f(y) > \alpha$ for every $y \in\gamma$.

There exists\vspace*{1pt} $z' \in A'$ such that $\rho(\vphi_D(z'),\eta/10)
\subset\vphi_D(A')$
and for every $\xi\in\rho(\vphi_D(z'),\eta/10)$,
every path from $\vphi_D^{-1}(\xi)$ to $\p D$ that
does not hit $\{ \zeta\in D \dvtx|\vphi_D(\zeta) - \xi| > \eps/
10 \}$ crosses $\gamma$
(as a continuous curve).
By the Koebe $1/4$ theorem 
and by the Koebe distortion theorem, 
there exist a finite set $Z \subset\C$ and $\eta'>0$, depending only
on $\eta$, such that for all
$\rho= \rho(\xi,\eta/10) \subset(1-\eta)\U$ and any
$D \in\D$, there exists $z \in Z$ with $\rho(z,\eta') \subset\vphi
_D^{-1}(\rho)$.
Thus, by\vadjust{\goodbreak} weak convergence and Proposition~\ref{prop: mu-contin
encoppasing a point},
for small enough $\delta_0$ (depending only on $\eta$),
there exists $z \in V_\delta(D)$ such that $\vphi_D(z) \in\rho
(\vphi_D(z'),\eta/10)$.
%
The probability that $S_z$ hits $\{ \zeta\in D \dvtx|\vphi_D(\zeta)
- \vphi_D(z)| > \eps/ 10 \}$
before exiting $D$ is at least $\min_{y \in\gamma} f(y) > \alpha$.
This is a contradiction to Proposition~\ref{prop: starting near boundary}.
\end{pf}
\begin{pf*}{Proof of Lemma~\ref{lem: starting near boundary}}
Let $\eta,\eta' > 0$ be small enough.
We show that if $\delta_0$ is small enough,
for every $D \in\D$, and for every $x \sim y \in V_\delta(D)$, we have
$|\vphi_D(x)-\vphi_D(y)| < \eta'$.

By the\vspace*{1pt} Koebe distortion theorem, 
using (\ref{eqn: der of vphi is bounded at 0}), for every $z \in
(1-\eta)\U$,
we have $|{\vphi_D^{-1\prime}}(z)| \geq c \eta$ for a constant $c > 0$.
By weak convergence, since $G$ is planar-irreducible, when $\delta_0$
tends to $0$,
the length of the edges of $G_\delta$ in $R \U$, for $R = \sup\{|z|
\dvtx z \in{\mathbf D}\}$, tends to $0$.
This implies that if $\delta_0$ is small enough,
for every $D \in\D$ and $y \sim x \in V_\delta(D)$ such that
$|\vphi_D(y)|,|\vphi_D(x)| \leq1 - \eta$,
we have $|\vphi_D(y) - \vphi_D(x)| \leq\eta'$.

It remains to consider $x$'s such that $|\vphi_D(x)| \geq1- \eta$.
As above, for small enough~$\delta_0$,
every $z \in[x,y]$ admits $|\vphi_D(z)| \geq1- 2\eta$.
Assume toward a contradiction that $|\vphi_D(x) - \vphi_D(y)| \geq
\eta'$.
Thus, by Proposition~\ref{prop: starting near boundary}
(using a similar argument to the one in Corollary~\ref{cor: starting
near boundary}),
there exists $\xi\in V_\delta(D)$ such that $1 - 4\eta\leq|\vphi
_D(\xi)| \leq1 - 2 \eta$
and the probability that $S_{\xi}$ hits the set $\{\zeta\in D \dvtx
|\vphi_D(\zeta) - \vphi_D(\xi)| > \eta'/3 \}$
before exiting $D$ is smaller than $1$.
However, since $G$ is planar-irreducible,
$S_{\xi}$ cannot cross $[x,y]$, so
the probability that $S_{\xi}$ hits the set $\{\zeta\in D \dvtx
|\vphi_D(\zeta) - \vphi_D(\xi)| > \eta'/3 \}$
before exiting $D$ is $1$, which is a contradiction.

The proof of the lemma follows by the strong Markov property,
and by applying Corollary~\ref{cor: starting near boundary} a finite
number of times.
\end{pf*}

\subsection{Exit probabilities are correct}

Let $D \in\D$. For $J \subset\p D$, denote by
$H(a,J;D)$ the probability that the natural random walk started at $a$
exits $D$ at~$J$; that is,
$H(a,J;D) = \sum_b H(a,b;D)$,
where the sum is over all $b \in\p V_{\delta}(D)$ such that $b \cap J
\neq\varnothing$.
\begin{lem} \label{lem: convergence from a}
For all $\eps, \alpha> 0$, for all $D \in\D$,
and for all $J = \vphi_D^{-1}(I)$ where $I \subset\p\U$ is an arc,
there exists $\delta_0 > 0$ such that for all $0 < \delta< \delta_0$
the following holds:

Let $a \in V_\delta(D)$ be such that $|\vphi_D(a)| \leq1-\eps$.
Then,
\[
\vert H(a,J;D) - \mathbb{P}_a [B(\tau) \in J ] \vert<
\alpha,
\]
where $B(\cdot)$ is a planar Brownian motion, and $\tau$ is the exit
time of $B(\cdot)$ from~$D$.
\end{lem}
\begin{pf}
Fix $\eps,\alpha,D$ and $J$ as above.
Denote $\vphi= \vphi_D$ and denote $\tau^{(a)} = \tau^{(a)}_D$.
Let $0 < \alpha_0 < 1$ be such that $\frac{(1+\alpha_0)^2}{(1-\alpha
_0)^2} = 1 + \frac{\alpha}{2}$.
Let $\eta> 0$ be small enough.
Denote $\mathcal{A} = \{ \frac{\eta}{4} (n + m \cdot i) \in
(1-\eps) \U\dvtx n,m \in\Z\}$.
The set $\mathcal{A}$ is finite, and there exists
${\tilde{a}} \in\mathcal{A}$ such that $\vphi(a) \in\rho({\tilde
{a}},\eta)$.
Denote $\rho= \vphi^{-1}(\rho(\tilde{a},\eta))$.\vadjust{\goodbreak}

We show that if $\eta,\delta_0$ is small enough, then
$\mathbb{P} [ S_x(\tau^{(x)}) \in J ] > (1-\alpha_0/2) \cdot\mathbb
{P} [ S_y(\tau^{(y)}) \in J ]$
for every $x,y \in V_\delta(\rho)$.
Define $h(z)$ to be the probability that $S_z(\tau^{(z)}) \in J$.
The map $h(\cdot)$ is harmonic in $V_\delta(D)$
with respect to the law of the natural random walk on $G_\delta$.
Since $h(\cdot)$ is harmonic,
there exists a path $\gamma$ from $y$ to $\p D$
such that $h(z) \geq h(y)$ for every $z \in\gamma$.
Since $h(\cdot)$ is nonnegative, harmonic and bounded,
\[
h(x) \geq\mathbb{P}\bigl[ S_x \bigl[ 0,\tau^{(x)}\bigr] \cap
\gamma\neq
\varnothing
\bigr] \cdot h(y).
\]
By Proposition~\ref{prop: close points intersect - fixed a},
since $G$ is planar,
$\mathbb{P}[ S_x [ 0,\tau^{(x)}] \cap\gamma\neq\varnothing] > 1 -
\alpha_0/2$
for small enough $\eta,\delta_0$.

Therefore, for small enough $\eta,\delta_0$,
%
%
\begin{equation} \label{eqn: a close to z}
\biggl\vert\frac{ \mathbb{P} [S_z(\tau^{(z)}) \in J ] }{ \mathbb{P}
[S_a(\tau^{(a)}) \in J ] } -1 \biggr\vert
< \alpha_0
\end{equation}
for every $z \in V_\delta(\rho)$.
In addition,
%
%
\begin{equation} \label{eqn: a close to z BM}
\biggl\vert\frac{ \mathbb{P}_z [B(\tau) \in J ] }{ \mathbb{P}_a
[B(\tau
) \in J ] } -1 \biggr\vert< \alpha_0
\end{equation}
for every $z \in\rho$.
By weak convergence and Proposition~\ref{prop: mu-contin encoppasing a point},
by the conformal invariance of Brownian motion,
we can choose $\delta_0$ so that
%
%
\begin{equation} \label{eqn: w.c. of 0+}
\biggl\vert\frac{ \mathbb{P} [\Theta_0(\rho) < \tau^{(0)} ,
S_0(\tau
^{(0)}) \in J ] }{ \mathbb{P}_0 [B[0,\tau] \cap\rho\neq\varnothing
, B(\tau) \in J ] } - 1 \biggr\vert< \alpha_0
\end{equation}
and
%
%
\begin{equation} \label{eqn: w.c. of 0}
\biggl\vert\frac{ \mathbb{P} [\Theta_0(\rho) < \tau^{(0)} ] }{
\mathbb
{P}_0 [B[0,\tau] \cap\rho\neq\varnothing] } - 1 \biggr\vert<
\alpha_0 .
\end{equation}

Combining (\ref{eqn: w.c. of 0+}) and (\ref{eqn: a close to z BM}),
\begin{eqnarray*}
&&\mathbb{P} \bigl[\Theta_0(\rho) < \tau^{(0)} , S_0\bigl(\tau
^{(0)}\bigr) \in J \bigr]\\
&&\qquad< (1+ \alpha_0) \mathbb{P}_0 \bigl[B[0,\tau] \cap\rho
\neq
\varnothing, B(\tau) \in J \bigr] \\
&&\qquad< (1+ \alpha_0)^2 \mathbb{P}_0 \bigl[B[0,\tau] \cap\rho
\neq
\varnothing\bigr] \mathbb{P}_a [B(\tau) \in J ]
\end{eqnarray*}
and combining (\ref{eqn: a close to z}) and (\ref{eqn: w.c. of 0}),
\begin{eqnarray*}
&&\mathbb{P} \bigl[\Theta_0(\rho) < \tau^{(0)} , S_0\bigl(\tau
^{(0)}\bigr) \in J \bigr]\\
&&\qquad> (1 - \alpha_0) \mathbb{P} \bigl[\Theta_0(\rho) < \tau
^{(0)}\bigr] \mathbb{P} \bigl[S_a\bigl(\tau^{(a)}\bigr) \in J
\bigr] \\
&&\qquad> (1 - \alpha_0)^2 \mathbb{P}_0 \bigl[B[0,\tau] \cap\rho
\neq
\varnothing\bigr] \mathbb{P} \bigl[S_a\bigl(\tau^{(a)}\bigr) \in
J \bigr].
\end{eqnarray*}
Thus, by the choice of $\alpha_0$,
\[
\mathbb{P} \bigl[S_a\bigl(\tau^{(a)}\bigr) \in J \bigr] <
(1+\alpha)
\mathbb{P}_a [B(\tau) \in J ].
\]
Similarly, since $1-\alpha< \frac{1}{1+\alpha/2}$,
\[
\mathbb{P} \bigl[S_a\bigl(\tau^{(a)}\bigr) \in J \bigr] >
(1-\alpha)
\mathbb{P}_a [B(\tau) \in J ].
\]
The lemma follows, since $\mathbb{P}_a [B(\tau) \in J ]
\leq1$.
\end{pf}

Using Lemma~\ref{lem: starting near boundary}, Lemma~\ref{lem:
convergence from a} yields the following.
\begin{lem} \label{lem: convergence from a for all D}
There exists a universal constant $c > 0$
such that for all $\alpha> 0$, there exists $\delta_0 > 0$
such that for all $0 < \delta< \delta_0$ the following holds:

Let $D \in\D$, and let $J = \vphi_D^{-1}(I)$, where $I \subset\p\U
$ is an arc of length at least $\alpha$.
Then,
$ H(0,J;D) \geq c \cdot\alpha$.
\end{lem}
\begin{pf}
Let $\eta> 0$ be small enough, and let $\tilde D$ be the $(1,\eta
)$-approximation of $D$
given by Proposition~\ref{prop: compact approx}.
Let $x \in\p\U$ be the center of $I$, and let $A = \rho(x,\alpha
/2) \cap\U$.
Let $\mathcal{I}$ be the finite family of arcs
of the form
$I = \{e^{i s} \dvtx\alpha j/8 \leq s \leq\alpha(j+1)/8
\}$
for $0 \leq j \leq16\pi/\alpha$.

Let $I' \in\mathcal{I}$ be so that $x \in I'$.
For every $\zeta\in I'$, since $|x-\zeta| \leq\alpha/8$
and since $|\vphi_D(\vphi^{-1}_{\tilde D}(\zeta)) - \zeta| \leq
\eta$,
we have $|x - \vphi_D(\vphi^{-1}_{\tilde D}(\zeta))| \leq\eta+
\alpha/8 < \alpha/4$
for $\eta< \alpha/8$.
Thus, $\operatorname{dist}(x,\vphi_D(\vphi^{-1}_{\tilde D}(I')) <
\alpha/4$.
As in the proof of Lemma~\ref{lem: starting near boundary},
if $\delta_0$ is small enough (independently of $D$),
for every $v \sim u \in V_\delta(D)$, we have $|\vphi_D(v) - \vphi
_D(u)| < \eta$.
Thus, by properties \hyperlink{item:com1}{(1)} and \hyperlink
{item:com3}{(3)} of
Proposition~\ref{prop: compact approx},
\[
H(0,J;D) \geq\mathbb{P} \bigl[S_0\bigl(\tau^{(0)}_{\tilde D}\bigr
) \in\vphi
^{-1}_{\tilde{D}}(I') \bigr]
\cdot\min_{y} \mathbb{P} \bigl[ S_y\bigl(\tau^{(y)}_{D}\bigr)
\in J \bigr],
\]
where the minimum is over
$y \in V_\delta(\rho(x,\alpha/2))$ such that $|\vphi_D(y)| \geq1-
2\eta$.
By weak convergence and Proposition~\ref{prop: mu-contin encoppasing a point},
if $\delta_0$ is small enough, we have that
$\mathbb{P} [S_0(\tau^{(0)}_{\tilde D}) \in\vphi^{-1}_{\tilde
{D}}(I') ]$
is at least a universal constant times $\alpha$.
By Lemma~\ref{lem: starting near boundary}, for small enough $\eta,
\delta_0$,
we have
$\min_{y} \mathbb{P} [ S_y(\tau^{(y)}_{D}) \in J ] \geq1/2$.
\end{pf}

\section{Convergence of Poisson kernel}
\label{sec: con to poison ker}

In this section we prove that one can approximate the discrete
Poisson kernel by the continuous Poisson kernel.

\subsection{\texorpdfstring{Proof of Lemma \protect\ref{lem: H(,) are close to lambda}}{Proof of Lemma 1.2}}

We begin with a proposition that is a ``special case'' of Lemma \ref
{lem: H(,) are close to lambda}
for a specific domain.
\begin{prop} \label{prop: H(,) are close to lambda}
Let $\eps, \alpha> 0$
and let $D\subsetneqq\C$ be a simply connected domain such that $0
\in D$.
Then, there exists $\delta_0$ such that for all $0 < \delta< \delta
_0$ the following holds:

Let $a \in V_\delta(D)$ be such that $\vert\vphi_D(a)\vert\leq
1-\eps
$, and let
$b \in\p V_{\delta}(D)$. Then,
\[
\biggl\vert\frac{ H^{(\delta)}(a,b;D)}{ H^{(\delta)}(0,b;D) } -
\lambda
(a,b;D) \biggr\vert\leq\alpha.
\]
\end{prop}

Roughly, Proposition~\ref{prop: H(,) are close to lambda} yields
Lemma~\ref{lem: H(,) are close to lambda}
by a compactness argument.\vadjust{\goodbreak}
\begin{pf*}{Proof of Lemma~\ref{lem: H(,) are close to lambda}}
Let $\alpha_1 > 0$ be small enough, and let $\tilde D$ be the
$(\eps,\alpha_1)$-approximation of $D$ given by Proposition \ref
{prop: compact approx}.
Let $\delta_0 > 0$ be small enough, and let $0 < \delta< \delta_0$.
Specifically, Proposition~\ref{prop: H(,) are close to lambda} holds
for $\tilde D$ with $\eps/2$ and
$\alpha_1$. Since $|\vphi_{\tilde{D}}(a)| \leq1- \eps/2$,
for every $\tilde b \in\p V_\delta(\tilde D)$,
\[
\biggl\vert\frac{H(a,\tilde b;\tilde D)}{H(0,\tilde b;\tilde D)} -
\lambda
(a,\tilde b;\tilde D)\biggr\vert\leq\alpha_1 .
\]
Since $\tilde D \subset D$, for every $x \in
V_\delta(\tilde D)$,
\[
H(x,b;D) = \sum_{\tilde b} H(x,\tilde b; \tilde D) \cdot H(\tilde b, b
; D) ,
\]
where the sum is over $\tilde b \in\p V_\delta(\tilde D)$,
and we abuse notation and use $H(\tilde b, b ; D)$ instead of
$H(\tilde{b}_2 , b ; D)$, where $\tilde b = (\tilde{b}_1,\tilde{b}_2)$
[for every $b' \in\p V_\delta(D)$, define $H(b',b;\break D) = {\mathbf1}_{\{
b =
b'\}}$].
Thus,
%
%
\begin{eqnarray}
\label{eqn: H lemma}
&&|H (a,b;D) - \lambda(a,b;D) \cdot H(0,b;D)| \nonumber\\
&&\qquad\leq
\sum_{\tilde b} H(\tilde b, b ; D) \cdot|H(a,\tilde b; \tilde D) -
\lambda(a,b;D) \cdot H(0,\tilde b; \tilde D) | \nonumber\\[-8pt]\\[-8pt]
&&\qquad\leq\sum_{\tilde b} H(\tilde b, b ; D) \cdot H( 0,\tilde b;
\tilde D) \cdot
|\lambda(a,\tilde b; \tilde D) - \lambda(a,b;D) |\nonumber\\
&&\qquad\quad{} + \alpha_1 \cdot
H(0,b;D).\nonumber
\end{eqnarray}
Let $\alpha_2,\alpha_3 > 0$ be small enough. Let $I \subset\p\U$
be an arc of length $\alpha_2$
centered at $\vphi_D(b)$. Denote $\tilde I = \vphi^{-1}_{\tilde D}(I)
\subset\p\tilde D$.
We use the following two claims.
\begin{clm}
\label{clm: point in tilde I}
For every $\tilde b \in\p V_{\delta}(\tilde D)$ such that $\tilde b
\cap\tilde I \neq\varnothing$,
$|\lambda(a,\tilde b; \tilde D) - \lambda(a,b;D)| \leq\alpha_3$.
\end{clm}
\begin{pf}
By the choice of $I$, $|\vphi_{\tilde D}(\tilde b) - \vphi_D(b)| \leq
\alpha_2$.
Since $a \in\tilde D$,
by property~\hyperlink{item:com4}{(4)} of Proposition~\ref{prop:
compact approx}
with $\xi= \vphi_{\tilde D}(a)$, we have
$|\vphi_{D}(a) - \vphi_{\tilde D}(a)| =
|\vphi_{D}(\vphi_{\tilde D}^{-1}(\xi)) - \xi| \leq\alpha_1$.
By the continuity\vspace*{1pt} of $\lambda(\cdot,\cdot;\U)$,
if $\alpha_1,\alpha_2$ are small enough,
$|\lambda(a,\tilde b; \tilde D) - \lambda(a,b;D)| \leq\alpha_3$.
\end{pf}
\begin{clm}
\label{clm: point out tilde I}
For every $\tilde b \in\p V_{\delta}(\tilde D)$ such that $\tilde b
\cap\tilde I = \varnothing$,
$H(\tilde b, b ; D) \leq\alpha_3 \cdot H(0,b;D)$.
\end{clm}
\begin{pf}
Assume that\vspace*{1pt} $\tilde b \notin\p V_\delta(D)$
[otherwise, $H(\tilde b,b;D) = 0$, since $\tilde b \cap I = \varnothing$].
In this case, $H(\tilde b, b ; D)$ is
$H({\tilde b}_2,b;D)$ where $\tilde{b}_2$ is the endpoint of $\tilde b$.
Denote $b' = \vphi_{\tilde D}(\tilde b) \in\p\U$. So $|b'-\vphi
_D(b)| \geq\alpha_2/10$.
By property \hyperlink{item:com4}{(4)} of Proposition~\ref{prop:
compact approx},
$|\vphi_D(\vphi^{-1}_{\tilde D}(b')) - b'| \leq\alpha_1$.
By weak convergence for small enough $\delta_0$,
the length of the edge $\tilde b$ is small enough.
Thus,\vspace*{1pt} by property \hyperlink{item:com5}{(5)} of
Proposition \ref
{prop: compact approx},
for small enough\vadjust{\goodbreak} $\delta_0$, we have $|\vphi_{D}(\tilde{b}_2) -
\vphi_D(\vphi^{-1}_{\tilde D}(b'))|
\leq\alpha_1$, which implies
$|\vphi_{D}(\tilde{b}_2) - b'| \leq2 \alpha_1$.
Therefore, $|\vphi_D(\tilde{b}_2)-\vphi_D(b)| \geq\alpha_2/10 - 2
\alpha_1 > \alpha_2 /20$,
for $\alpha_1 < \alpha_2/40$.\vspace*{1pt}

Denote $\xi= \vphi_D(\tilde{b}_2)$, and $A = \{ x \in\U\dvtx
|x-\xi| > \alpha_2/50\}$.
Also denote $M = \max_y H(y,b;D)$,
where the maximum is over $y \in V_\delta(D)$ such that
$|\vphi_D(y) - \vphi_D(b)| \geq\alpha_2/50$.
As in the proof of Lemma~\ref{lem: starting near boundary},
if $\delta_0$ is small enough,
for every $v \sim u \in V_\delta(D)$, we have $|\vphi_D(v) - \vphi
_D(u)| < \alpha_2/100$.
Thus, $H(\tilde b,b;D)$ is at most
$M$ times the probability that $\vphi_D \circ S_{\tilde{b}_2}$ hits $A$.

Since $|\xi- \vphi_D(\vphi^{-1}_{\tilde D}(b'))|
\leq\alpha_1$, using property \hyperlink{item:com3}{(3)} of
Proposition \ref
{prop: compact approx},
$|\xi| \geq1- 2\alpha_1$.
Let $\alpha_4 > 0$ be small enough.
Using Lemma~\ref{lem: starting near boundary},
for $\alpha_1$ small enough,
the probability that $\vphi_D \circ S_{\tilde{b}_2}$ hits $A$ is at most
$\alpha_4$.\vspace*{2pt}

We show that $M \leq C \cdot H(0,b;D)$, for some $C= C(\alpha_2) > 0$.
Let $y \in V_\delta(D)$ be such that $|\vphi_D(y) - \vphi_D(b)| \geq
\alpha_2/50$.
The map $H(\cdot,b;D)$ is harmonic with respect to the law
of the natural random walk on $G_\delta$.
Thus, there exists a path $\gamma$ from $y$ to $b$ in $V_\delta(D)$
such that
for every $z \in\gamma$, $H(z,b;D) \geq H(y,b;D)$.
Since $H(\cdot,b;D)$ is nonnegative, harmonic and bounded,
\[
H(0,b;D) \geq\mathbb{P} \bigl[ S_0\bigl[0,\tau^{(0)}_D\bigr] \cap
\gamma\neq
\varnothing\bigr] \cdot H(y,b;D) .
\]
Therefore, we need to show that
$p = \mathbb{P} [ S_0[0,\tau^{(0)}_D] \cap\gamma\neq\varnothing]$
can be bounded from below
by a function of $\alpha_2$.

Think of $\gamma$ as a continuous curve, and denote
$\gamma' = \{ \zeta\in\gamma\dvtx|\vphi_D(\zeta) - \vphi_D(b)|
\leq\alpha_2 / 50 \}$.
Denote $D ' = D \setminus\gamma'$.
By the conformal invariance of the harmonic measure,
the length of the arc $\vphi_{D'}(\gamma')$ is at least a universal
constant times~$\alpha_2$.
Also, for small enough $\alpha_2$, we have $\operatorname{rad}(D')
\geq1/4$.
Thus, by Lemma~\ref{lem: convergence from a for all D} applied to $D'$
(using Lemma~\ref{lem: convergence from a for all D} with $\D' =
\{2D \dvtx D \in\D\}$),
$p$ is at least a universal constant times~$\alpha_2$. Set $C(\alpha
_2) = \frac{1}{p}$.

Setting $\alpha_4 \cdot C(\alpha_2) \leq\alpha_3$, the proof is complete.
\end{pf}

By Claims~\ref{clm: point in tilde I} and~\ref{clm: point out tilde I},
\begin{eqnarray*}
(\ref{eqn: H lemma}) & \leq&\sum_{\tilde b \dvtx\tilde b \cap
\tilde
I \neq\varnothing} H(\tilde b, b ; D) \cdot H( 0,\tilde b; \tilde D)
\cdot
|\lambda(a,\tilde b; \tilde D) - \lambda(a,b;D) | \\
&&{} + \sum_{\tilde b \dvtx\tilde b \cap\tilde I = \varnothing}
H(\tilde b, b ; D) \cdot H( 0,\tilde b; \tilde D) \cdot
|\lambda(a,\tilde b; \tilde D) - \lambda(a,b;D) |\\
&&{} + \alpha_1 \cdot
H(0,b;D) \\
& \leq&\alpha_3 \cdot\sum_{\tilde b} H(\tilde b, b ; D) \cdot H(
0,\tilde b; \tilde D)
+ \alpha_3 \cdot H(0,b;D) \cdot\sum_{\tilde b} H( 0,\tilde b; \tilde
D)\\
&&{} + \alpha_1 \cdot H(0,b;D) \\
& \leq&( 2 \alpha_3 + \alpha_1 ) \cdot H(0,b;D).
\end{eqnarray*}
Choosing $2\alpha_3+\alpha_1 < \alpha$ completes the proof.\vadjust{\goodbreak}
\end{pf*}

\subsection{\texorpdfstring{Proof of Proposition \protect\ref{prop: H(,) are close to lambda}}{Proof of Proposition 5.1}}

Fix $\eps, \alpha> 0$.
Let $N$ be a large enough integer so that
%
%
\begin{equation}
\label{eqn: cond on N} \bigl(1-c_1(\eps,\alpha)\bigr)^N < \frac
{\alpha}{8 c_2},
\end{equation}
where $c_1(\eps,\alpha) > 0$ is given below in Proposition \ref
{prop: a inter 0},
and $c_2 > 0$ is the universal constant given below in Proposition \ref
{prop: a exits b, like w}.
Let $\beta> 0$ be small enough so that
%
%
\begin{equation}
\label{eqn: cond on beta}
\beta< \frac{\eps}{50 \pi K 5^N} \quad\mbox{and}\quad
\beta< \frac{\alpha\eps^5}{16 c_3},
\end{equation}
where $K>5$ is the universal constant from
Lemmas~\ref{lem: BM from 0 several} and~\ref{lem: BM from 0 hits I+ I-},
and $c_3 > 0$ is the universal constant given below in (\ref{eqn: def
of c3}),
and let $r = 2 \pi\beta$.
Let $\eta> 0$ be given by Proposition~\ref{prop: close points
intersect - fixed a} with
$\alpha$ equals $\frac{\alpha\eps^2}{16}$.
Let $\delta_0 > 0$ be small enough
(to be determined below), and let $0 < \delta< \delta_0$.

Denote $\vphi= \vphi_D$.
Denote
$\mathcal{A} = \{ \frac{\eps}{100} (n + m \cdot i) \in(1-\eps
) \U\dvtx n,m \in\Z\}$.
The set $\mathcal{A}$ is finite,
and there exists ${\tilde{a}} \in\mathcal{A}$ such that $\vphi(a)
\in\rho({\tilde{a}},\eps/40)$.
Denote
$\mathcal{B} = \{ e^{\pi\beta n i /20} \dvtx0 \leq n \leq100/\beta
\}$.
The set $\mathcal{B}$ is finite, and there exists $\tilde b \in
\mathcal{B}$ such that
$|\vphi(b) - \tilde b| \leq\beta/10$.
Denote
$I = \{ {\tilde b} \cdot e^{i t} \dvtx- \pi\beta\leq t \leq\pi
\beta\}$,
and denote $J = \vphi^{-1} (I)$.
Roughly, $b$ is an edge in the middle of the small interval $J$.

For $j=1,2, \ldots,N$, let $R_j = 5^j K r$,
let $\xi_j = {\tilde b} (1-3R_j)$,
and let $\rho_j = \rho(\xi_j, \eta^3 R_j)$.
For $z \in V_\delta({D})$, define
\[
T^{(z)}_j = \min\{ t \geq0 \dvtx| \vphi(S_z(t))-{\tilde b}
| \leq R_j \}.
\]
On the event $\{S_z(\tau^{(z)}) \in J \}$, we have
$T^{(z)}_N \leq T^{(z)}_{N-1} \leq\cdots\leq T^{(z)}_1 \leq\tau^{(z)}$.
Let $E_j^{(z)}$ be the event
\[
E_j^{(z)} = \bigl\{ \vphi\circ S_0\bigl[T^{(0)}_{j+1},T^{(0)}_j\bigr
] \cap\rho
_j \neq\varnothing\bigr\} \cap
\bigl\{ \vphi\circ S_z\bigl[T^{(z)}_{j+1}, T^{(z)}_j\bigr] \cap
\rho_j \neq
\varnothing\bigr\}.
\]
Denote $E_j = E_j^{(a)}$.\vspace*{2pt}

We use the following three propositions.
\begin{prop} \label{prop: a inter 0}
Let $1 \leq j \leq N$. Then,
\[
\mathbb{P} \bigl[E_j \mid\ov{E}_{j+1} , \ldots, \ov{E}_{N} ,
S_0\bigl(\tau
^{(0)}\bigr) \in J , S_a\bigl(\tau^{(a)}\bigr) \in J \bigr] \geq c_1
\]
for $c_1 = c_1(\eps,\alpha) > 0$.
\end{prop}
\begin{prop} \label{prop: a exits b, like w}
There exists a universal constant $c_2 > 0$ such that
for every $z \in\{0,a\}$,
\begin{eqnarray*}
&&\mathbb{P} \bigl[ S_z\bigl(\tau^{(z)}\bigr) = b \mid\ov{E}_1 ,
\ldots,
\ov{E}_{N} , S_0\bigl(\tau^{(0)}\bigr) \in J,
S_a\bigl(\tau^{(a)}\bigr) \in J \bigr] \\
&&\qquad\leq c_2 \cdot
\mathbb{P} \bigl[S_0\bigl(\tau^{(0)}\bigr)
= b \mid S_0\bigl(\tau^{(0)}\bigr) \in J \bigr] .
\end{eqnarray*}
\end{prop}
\begin{prop}
\label{prop: cond on hiting rhoj}
For every $j = 1,\ldots,N$,
\[
\biggl\vert\frac{ \mathbb{P} [S_a(\tau^{(a)}) = b \mid\exit,
E_j, \ov
{E}_{j+1},\ldots,\ov{E}_{N} ] } {\mathbb{P} [S_0(\tau^{(0)}) = b
\mid\exit, E_j, \ov{E}_{j+1},\ldots,\ov{E}_{N} ] } - 1 \biggr
\vert
\leq\frac{\alpha\eps^2}{4} .
\]
\end{prop}

Before proving the three propositions above,
we show how they imply Proposition~\ref{prop: H(,) are close to lambda}.
Let $z \in\{0,a\}$.
Write
%
%
\begin{eqnarray}
\label{eqn: decompose H}
H(z,b) &=& H^{(\delta)}(z,b;D) \nonumber\\[-8pt]\\[-8pt]
&=& \mathbb{P} \bigl[S_z\bigl(\tau^{(z)}\bigr) = b \mid
S_z\bigl(\tau^{(z)}\bigr) \in J \bigr]
\cdot\mathbb{P} \bigl[S_z\bigl(\tau^{(z)}\bigr) \in J \bigr]
.\nonumber
\end{eqnarray}
By Lemma~\ref{lem: convergence from a},
by (\ref{eqn: harmonic measure}),
and since $|\vphi(z)| \leq1 - \eps$,
%
%
\begin{equation}
\label{eqn: def of c3}
\bigl| \mathbb{P} \bigl[S_z\bigl(\tau^{(z)}\bigr) \in J\bigr] -
\lambda
(z,b) \cdot\beta\bigr| \leq c_3 \frac{\beta^2}{ \eps^3}
\end{equation}
for a universal constant $c_3 > 0$, which implies
%
%
\begin{equation}
\label{eqn: a divide 0}
\biggl\vert\frac{\mathbb{P} [S_a(\tau^{(a)}) \in J]} {\mathbb{P}
[S_0(\tau^{(0)}) \in J]} - \lambda(a,b) \biggr\vert<
\frac{\alpha}{4} .
\end{equation}
Denote
$\exit= \{S_0(\tau^{(0)}) \in J\} \cap\{S_a(\tau
^{(a)}) \in J\}$,
and denote
$\inter= E_1 \cup E_2 \cup\cdots\cup E_N$.
Since
\begin{eqnarray*}
&&\mathbb{P} \bigl[S_z\bigl(\tau^{(z)}\bigr) = b \mid S_z\bigl
(\tau^{(z)}\bigr) \in J \bigr] \\
&&\qquad= \sum_{j = 1}^N \mathbb{P} \bigl[S_z\bigl(\tau
^{(z)}\bigr) = b \mid\exit, E_j,
\ov{E}_{j+1},\ldots,\ov{E}_{N} \bigr] \\
&&\qquad\quad\hspace*{13.7pt}{}\times\mathbb{P} [E_j, \ov
{E}_{j+1},\ldots,\ov{E}_{N} \mid\exit]
\\
&&\qquad\quad{} + \bigl[S_z\bigl(\tau^{(z)}\bigr) = b \mid\exit
, \ninter\bigr] \cdot
\mathbb{P} [\ninter\mid\exit] ,
\end{eqnarray*}
we have
%
%
\begin{eqnarray}
\label{eqn: a compared to 0 first}
&&\bigl| \mathbb{P} \bigl[ S_a\bigl(\tau^{(a)}\bigr) = b \mid
S_a\bigl(\tau
^{(a)}\bigr) \in J \bigr] -
\mathbb{P} \bigl[S_0\bigl(\tau^{(0)}\bigr) = b \mid S_0\bigl(\tau
^{(0)}\bigr) \in J \bigr]
\bigr| \nonumber\\
&&\qquad\leq\sum_{j = 1}^N \mathbb{P} [E_j, \ov{E}_{j+1},\ldots
,\ov
{E}_{N} \mid\exit] \nonumber\\
&&\qquad\quad\hspace*{13pt}{} \times\bigl| \mathbb{P} \bigl
[S_a\bigl(\tau^{(a)}\bigr) = b \mid\exit, E_j, \ov
{E}_{j+1},\ldots,\ov{E}_{N} \bigr] \\
&&\qquad\quad\hspace*{26pt}{} - \mathbb{P} \bigl[S_0\bigl(\tau
^{(0)}\bigr) = b \mid\exit, E_j, \ov
{E}_{j+1},\ldots,\ov{E}_{N} \bigr] \bigr|
\nonumber\\
&&\qquad\quad{} + 2 \max_{z \in\{0,a\}} \bigl[S_z\bigl(\tau
^{(z)}\bigr) = b \mid\exit, \ninter
\bigr] \cdot\mathbb{P} [\ninter\mid\exit] .\nonumber
\end{eqnarray}
By Propositions~\ref{prop: a inter 0} and~\ref{prop: a exits b, like w},
%
%
\begin{eqnarray}
\label{eqn: 0 and a to 0}
&&2 \max_{z \in\{0,a\}} \bigl[S_z\bigl(\tau^{(z)}\bigr) = b \mid
\exit, \ninter\bigr]
\cdot\mathbb{P} [\ninter\mid\exit]
\nonumber\\[-8pt]\\[-8pt]
&&\qquad< \mathbb{P} \bigl[S_0\bigl(\tau^{(0)}\bigr) = b \mid
S_0\bigl(\tau
^{(0)}\bigr) \in J \bigr] \cdot\frac{\alpha}{4}.\nonumber
\end{eqnarray}

Plugging Proposition~\ref{prop: cond on hiting rhoj} and (\ref{eqn: 0
and a to 0}) into (\ref{eqn: a compared to 0 first}),
\begin{eqnarray*}
&&
\bigl| \mathbb{P} \bigl[ S_a\bigl(\tau^{(a)}\bigr) = b \mid
S_a\bigl(\tau
^{(a)}\bigr) \in J \bigr] -
\mathbb{P} \bigl[S_0\bigl(\tau^{(0)}\bigr) = b \mid S_0\bigl(\tau
^{(0)}\bigr) \in J \bigr]
\bigr| \\
&&\qquad< \frac{\alpha\eps^2 }{2} \cdot\mathbb{P} \bigl[S_0\bigl
(\tau^{(0)}\bigr) = b
\mid S_0\bigl(\tau^{(0)}\bigr) \in J \bigr] .
\end{eqnarray*}
Thus, plugging (\ref{eqn: a divide 0}) into (\ref{eqn: decompose H}),
\[
\biggl\vert\frac{H(a,b)}{H(0,b)} - \lambda(a,b)\biggr\vert<
\frac{\alpha}{4} \biggl( 1 + \frac{\alpha\eps^2}{2} \biggr) +
\frac
{\alpha\eps^2}{2} \lambda(a,b) < \alpha.
\]
\subsection{\texorpdfstring{Proof of Proposition \protect\ref{prop: a inter 0}}{Proof of Proposition 5.4}}

For the rest of this proof denote by $E^{(z)}$ the event
\[
E^{(z)} = \ov{E}_{j+1}^{(z)} \cap\cdots\cap\ov{E}_{N}^{(z)} \cap
\bigl\{S_0\bigl(\tau^{(0)}\bigr) \in J\bigr\} \cap
\bigl\{S_z\bigl(\tau^{(z)}\bigr) \in J\bigr\},
\]
and denote $E = E^{(a)}$.
We show that
%
%
\begin{equation}
\label{eqn: second Sa}
\mathbb{P} \bigl[\vphi\circ S_a\bigl[T^{(a)}_{j+1},T^{(a)}_j\bigr]
\cap\rho_j \neq
\varnothing\mid E \bigr]
\end{equation}
is at least a constant (that may depend on $\eps$ and $\alpha$).
This implies the proposition, since $S_0$ and $S_a$ are independent
(and since the same argument holds for $0$ as well).
\begin{clm}
\label{clm: cut around a for first}
There exists a set of vertices $U \subset V_{\delta}({D})$ such that:
\begin{itemize}
\item Every path from $\vphi^{-1}(\rho(\tilde a,\eps/40))$
to the boundary of ${D} \setminus\vphi^{-1}((1- \eps/2) \U)$
in $G_{\delta}$ goes through $U$.
\item For every $u \in U$, we have
$\mathbb{P} [\vphi\circ S_u[T^{(u)}_{j+1},T^{(u)}_j] \cap\rho_j
\neq\varnothing\mid E^{(u)} ] \geq c_1$ with $c_1 = c_1(\eps,\alpha
) > 0$.
\end{itemize}
\end{clm}
\begin{pf}
Assume toward a contradiction that such a set does not exist.
Since $G$ is planar-irreducible,
there exists a path $Y$ from $\vphi^{-1}(\rho(\tilde a,\eps/40))$ to
the boundary of
$D \setminus\vphi^{-1}((1- \eps/2) \U)$ such that
for every vertex $y$ in $Y$,
%
%
\begin{equation}
\label{eqn: atc y}
\mathbb{P} \bigl[\vphi\circ S_y\bigl[T^{(y)}_{j+1},T^{(y)}_j\bigr]
\cap\rho_j \neq
\varnothing\mid E^{(y)} \bigr] < c_1(\eps,\alpha) .
\end{equation}

Define an auxiliary random walk $L$;
let $L(\cdot)$ be a natural random walk started at~$0$
(independent of $S_0$), and let $\tau^{(L)}$ be the exit time of
$L(\cdot)$ from ${D}$.
For \mbox{$j \leq k \leq N$}, let
\[
T_k^{(L)} = \min\{ t \geq0 \dvtx| \vphi(L(t))-\tilde{b} | \leq
R_k \},
\]
let
\[
E_k^{(L)} = \bigl\{ \vphi\circ S_0\bigl[T^{(0)}_{k+1},T^{(0)}_k\bigr
] \cap\rho_k
\neq\varnothing\bigr\} \cap
\bigl\{ \vphi\circ L\bigl[T^{(L)}_{k+1}, T^{(L)}_k\bigr] \cap\rho
_k \neq
\varnothing
\bigr\}
\]
and let
\[
E^{(L)} = \ov{E}_{j+1}^{(L)} \cap\cdots\cap\ov{E}_{N}^{(L)} \cap
\bigl\{S_0\bigl(\tau^{(0)}\bigr) \in J\bigr\} \cap\bigl\{L\bigl
(\tau^{(L)}\bigr) \in
J\bigr\}.
\]

Consider
%
%
\begin{equation}
\label{eqn: atc 0}
\mathbb{P} \bigl[L\bigl[0,T^{(L)}_{N}\bigr] \cap Y \neq\varnothing
, \vphi\circ
L\bigl[T^{(L)}_{j+1},T^{(L)}_j\bigr] \cap\rho_j \neq\varnothing
\mid E^{(L)}
\bigr] .
\end{equation}
By (\ref{eqn: atc y}), and by the strong Markov property,
we have $\mbox{(\ref{eqn: atc 0})} < c_1(\eps,\alpha)$.
On the other hand,
by weak convergence and Proposition~\ref{prop: mu-contin encoppasing a point},
by Lemma~\ref{lem: CONDITIONAL prob. of intersection a ball},
by Proposition~\ref{prop: encompassing a point},
and by the planarity of $G$,
\[
\mbox{(\ref{eqn: atc 0})} \geq
\mathbb{P} \bigl[\vphi\circ L[0,T'] \around^{(\eps/2)} \tilde a ,
\vphi
\circ L\bigl[T^{(L)}_{j+1},T^{(L)}_j\bigr] \cap\rho_j \neq
\varnothing\mid
E^{(L)} \bigr]
\geq c_2,
\]
where $T'$ is the first time $L(\cdot)$ hits
the set $\{z \in D \dvtx\vert\vphi(z)\vert\geq1-\eps/2\}$, and
$c_2 =
c_2(\eps,\alpha) > 0$.
This is a contradiction for $c_1 = c_2$.
\end{pf}

By Claim~\ref{clm: cut around a for first}, and by the strong Markov property,
(\ref{eqn: second Sa}) is a convex combination of
\[
\mathbb{P} \bigl[\vphi\circ S_u\bigl[T^{(u)}_{j+1},T^{(u)}_j\bigr]
\cap\rho_j \neq
\varnothing\mid E^{(u)} \bigr] \qquad\mbox{for $u \in U$} ,
\]
which implies that $\mbox{(\ref{eqn: second Sa})} \geq c_1(\eps
,\alpha)$.

\subsection{\texorpdfstring{Proof of Proposition \protect\ref{prop: a exits b, like w}}{Proof of Proposition 5.5}}

We use the following lemma, which is a variant of Harnack's inequality.
\begin{lem} \label{lem: Harnack for Pr[x -- b | I ]}
There exists a universal constant $c>0$ such that the following holds:

Let $w \in V_{\delta}(D)$ be such that $|\vphi(w)- \tilde b| \geq K r$.
If $\mathbb{P} [S_w(\tau^{(w)}) \in J ] > 0$, then
\[
\mathbb{P} \bigl[S_w\bigl(\tau^{(w)}\bigr) = b \mid S_w\bigl(\tau
^{(w)}\bigr) \in J \bigr] \leq
c \cdot\mathbb{P} \bigl[S_0\bigl(\tau^{(0)}\bigr) = b \mid
S_0\bigl(\tau^{(0)}\bigr) \in J
\bigr] .
\]
\end{lem}

Before proving the lemma, we show how the lemma implies
Proposition~\ref{prop: a exits b, like w}.
\begin{pf*}{Proof of Proposition~\ref{prop: a exits b, like w}}
Denote by $W$ the set of $w \in V_{\delta}(D)$ such that $|\vphi(w)-
\tilde b| \geq K r$
and $\mathbb{P} [S_w(\tau^{(w)}) \in J ] > 0$.
As in the proof of Lemma~\ref{lem: starting near boundary},
if $\delta_0$ is small enough,
for every $v \sim u \in V_\delta(D)$, we have $|\vphi(v) - \vphi(u)|
< \beta$.
By the strong Markov property,
\[
\mathbb{P} \bigl[ S_z\bigl(\tau^{(z)}\bigr) = b \mid\ov{E}_1 ,
\ldots
, \ov{E}_{N} , S_0\bigl(\tau^{(0)}\bigr) \in J , S_a\bigl(\tau
^{(a)}\bigr) \in
J \bigr]
\]
is at most
%
%
\begin{equation}
\label{eqn: max on W}
\max_{w \in W} \mathbb{P} \bigl[S_w\bigl(\tau^{(w)}\bigr) = b
\mid S_w\bigl(\tau
^{(w)}\bigr) \in J \bigr] .
\end{equation}
Lemma~\ref{lem: Harnack for Pr[x -- b | I ]} implies the proposition.\vadjust{\goodbreak}
\end{pf*}
\begin{pf*}{Proof of Lemma~\ref{lem: Harnack for Pr[x -- b | I ]}}
Let
\[
I_{+} = \{ \tilde b \cdot e^{it} \dvtx\pi\beta/2 \leq t \leq\pi
\beta\}
\quad\mbox{and}\quad
I_{-} = \{ \tilde b \cdot e^{it} \dvtx- \pi\beta\leq t \leq- \pi
\beta/2 \} .
\]
Let $J_+ = \vphi^{-1}(I_+)$ and $J_- = \vphi^{-1}(I_-)$.
Let
$U = \{ x \in{D} \dvtx|\vphi(x)- \tilde b| \geq Kr \}$,
let $\xi= \tilde b \cdot(1-3r)$, and let $\rho= \rho(\xi,r/20)$.

We use the following claim and its corollary.
\begin{clm}
\label{clm: exist good zetas} There exists a universal constant $c_1 >
0$ such that the following holds:

\begin{longlist}[(2)]
\item[(1)]\hypertarget{item:good0} There exists $x_0 \in V_\delta({
D}) \cap
\vphi^{-1}(\rho)$ such that
\[
\mathbb{P} \bigl[\vphi\circ S_{x_0}\bigl[0,\tau^{(x_0)}\bigr]
\around^{(r)} \xi,
S_{x_0}\bigl[0,\tau^{(x_0)}\bigr] \cap U = \varnothing\mid
S_{x_{0}}\bigl(\tau
^{(x_{0})}\bigr) \in J \bigr] \geq c_1 .
\]
\item[(2)]\hypertarget{item:good+} There exists $x_{+} \in V_\delta
({D}) \cap
\vphi^{-1}(\rho)$ such that
\[
\mathbb{P} \bigl[S_{x_{+}}\bigl(\tau^{(x_{+})}\bigr) \in J_{+} ,
S_{x_{+}}\bigl[0,\tau
^{(x_{+})}\bigr] \cap U = \varnothing\mid S_{x_{+}}\bigl(\tau
^{(x_{+})}\bigr) \in
J \bigr] \geq c_1.
\]
\item[(3)]\hypertarget{item:good-} There exists $x_{-} \in V_\delta
({D}) \cap
\vphi^{-1}(\rho)$ such that
\[
\mathbb{P} \bigl[S_{x_{-}}\bigl(\tau^{(x_{-})}\bigr) \in J_{-} ,
S_{x_{-}}\bigl[0,\tau
^{(x_{-})}\bigr] \cap U = \varnothing\mid S_{x_{-}}\bigl(\tau
^{(x_{-})}\bigr) \in
J \bigr] \geq c_1.
\]
\end{longlist}
\end{clm}
\begin{pf}
We first prove \hyperlink{item:good0}{(1)}.
Consider
%
%
\begin{eqnarray}
\label{eqn: atac for 0}
&&\mathbb{P} \bigl[\vphi\circ S_{0}\bigl[\T_0(\vphi^{-1}(\rho
)),\tau^{(0)}\bigr]
\around^{(r)} \xi,\nonumber\\[-8pt]\\[-8pt]
&&\qquad\hspace*{0pt} S_{0}\bigl[\T_0(\vphi^{-1}(\rho)),\tau
^{(0)}\bigr] \cap
U = \varnothing\mid S_0\bigl(\tau^{(0)}\bigr) \in J \bigr] .\nonumber
\end{eqnarray}
%
First, by weak convergence and Proposition~\ref{prop: mu-contin
encoppasing a point},
using Lemma~\ref{lem: BM from 0 several}, we have $\mbox{(\ref{eqn: atac
for 0})} \geq c_1$,
for a universal constant $c_1 > 0$.
Second, by the strong Markov property,
\[
\mbox{(\ref{eqn: atac for 0})}
\leq\max_{x} \mathbb{P} \bigl[\vphi\circ S_x\bigl[0,\tau
^{(x)}\bigr] \around
^{(r)} \xi, S_{x}\bigl[0,\tau^{(x)}\bigr] \cap U = \varnothing\mid
S_x\bigl(\tau
^{(x)}\bigr) \in J \bigr] ,
\]
where the maximum is over $x$ in $V_\delta({D}) \cap\vphi^{-1}(\rho
)$ such that \mbox{$\mathbb{P} [ S_x(\tau^{(x)}) \in J] > 0$}.

For the proof of property \hyperlink{item:good+}{(2)} we consider
$\{S_x(\tau^{(x)}) \in J_{+}\}$ instead of $\{\vphi
\circ S_x[0,\tau^{(x)}] \around^{(r)} \xi\}$, and use the
same argument
with Lemma~\ref{lem: BM from 0 hits I+ I-}.
Similarly, for property \hyperlink{item:good-}{(3)} we consider $\{
S_x(\tau^{(x)}) \in J_{-}\}$.
\end{pf}
\begin{cor}
\label{cor: zeta 0 exist J+- good}
There exists a universal constant $c_2 > 0$ such that the following holds:

There exists $x_0 \in V_\delta({D}) \cap\vphi^{-1}(\rho)$ such that
\[
\mathbb{P} \bigl[S_{x_{0}}\bigl(\tau^{(x_{0})}\bigr) \in J_{+} ,
S_{x_{0}}\bigl[0,\tau
^{(x_{0})}\bigr] \cap U = \varnothing\mid S_{x_{0}}\bigl(\tau
^{(x_{0})}\bigr) \in
J \bigr] \geq c_2
\]
and
\[
\mathbb{P} \bigl[S_{x_{0}}\bigl(\tau^{(x_{0})}\bigr) \in J_{-} ,
S_{x_{0}}\bigl[0,\tau
^{(x_{0})}\bigr] \cap U = \varnothing\mid S_{x_{0}}\bigl(\tau
^{(x_{0})}\bigr) \in
J \bigr] \geq c_2.
\]
\end{cor}
\begin{pf}
Let $x_0,x_{+},x_{-}$ be as given in Claim~\ref{clm: exist good zetas}.
We prove the first inequality for $x_0$,
the proof of the second one is similar.
Define
\[
h(z) =
\mathbb{P} \bigl[S_{z}\bigl(\tau^{(z)}\bigr) \in J_{+} ,
S_{z}\bigl[0,\tau^{(z)}\bigr] \cap
U = \varnothing\mid S_{z}\bigl(\tau^{(z)}\bigr) \in J \bigr] .
\]
The map $h(\cdot)$ is harmonic,
and so there exists a path $\gamma$ from $x_{+}$ to $\p D$ such that
$h(z) \geq h(x_{+})$ for every $z \in\gamma$.
Since $h(\cdot)$ is nonnegative, harmonic and bounded,
by Claim~\ref{clm: exist good zetas},
\begin{eqnarray*}
h(x_0) & \geq&\mathbb{P} \bigl[S_{x_0}\bigl[0,\tau^{(x_0)}_{D
\setminus U} \bigr]
\cap\gamma\neq\varnothing\mid S_{x_0}\bigl(\tau^{(x_0)}\bigr)
\in J \bigr]
\cdot h(x_{+}) \\
& \geq&\mathbb{P} \bigl[\vphi\circ S_{x_0}\bigl[0,\tau
^{(x_0)}\bigr] \around
^{(r)} \xi,\\
&&\hspace*{9.9pt} S_{x_0}\bigl[0,\tau^{(x_0)}\bigr] \cap U =
\varnothing\mid
S_{x_{0}}\bigl(\tau^{(x_{0})}\bigr) \in J \bigr]
\cdot h(x_{+})
\\
& \geq& c_2 .
\end{eqnarray*}
%
%
%
\upqed
\end{pf}

Back to the proof of Lemma~\ref{lem: Harnack for Pr[x -- b | I ]}.
For $y \in V_\delta({D})$, define
\[
p(y) = \cases{
\mathbb{P} \bigl[S_y\bigl(\tau^{(y)}\bigr) = b \mid S_y\bigl(\tau
^{(y)}\bigr) \in J \bigr] ,
&\quad
if $\mathbb{P} \bigl[S_y\bigl(\tau^{(y)}\bigr) \in J \bigr] > 0$,
\cr
0, &\quad otherwise.}
\]
Since $p(\cdot)$\vspace*{1pt} is harmonic, 
there exists a path $\gamma$ from $w$ to $b$ such that $p(z) \geq
p(w)$ for every $z \in\gamma$.
Let $x_0$ be the vertex given by Corollary~\ref{cor: zeta 0 exist J+- good}.
By the choice of $\tilde b$, $\vphi(b) \in I$ and $\vphi(b) \notin
I_+ \cup I_-$.
Thus, since $w \in U$, assume without loss of generality that every
path from $x_0$ to $J_{+}$
that does not intersect $U$ crosses $\gamma$ (otherwise, this holds
for $J_{-}$).
Thus, since $p(\cdot)$ is nonnegative, harmonic and bounded,
by Corollary~\ref{cor: zeta 0 exist J+- good},
%
%
\begin{eqnarray}
\label{eqn: x0 at least 0}
p(x_0) & \geq&\mathbb{P} \bigl[S_{x_0}\bigl[0,\tau^{(x_0)}\bigr]
\cap\gamma\neq
\varnothing\mid S_{x_0}\bigl(\tau^{(x_0)}\bigr) \in J \bigr] \cdot
p(w) \nonumber\\
& \geq&\mathbb{P} \bigl[S_{x_{0}}\bigl(\tau^{(x_{0})}\bigr) \in J_{+},
S_{x_{0}}\bigl[0,\tau^{(x_{0})}\bigr] \cap U = \varnothing\mid
S_{x_{0}}\bigl(\tau
^{(x_{0})}\bigr) \in J \bigr] \cdot p(w)
\hspace*{-28pt}\\
& \geq& c_2 \cdot p(w) ,\nonumber
\end{eqnarray}
where $c_2 > 0$ is a constant.

Similarly, there exists a path $\gamma$ from $x_0$ to $b$
(we abuse notation and use $\gamma$ again)
such that $p(z) \geq p(x_0)$ for every $z \in\gamma$.
Since $G$ is planar-irreducible,
every path from $0$ that encompasses $\vphi^{-1}(\rho)$ crosses
$\gamma$.
Since $p(\cdot)$ is nonnegative, harmonic and bounded,
\begin{eqnarray*}
p(0) & \geq&\mathbb{P} \bigl[S_{0}\bigl[0,\tau^{(0)}\bigr] \cap
\gamma\neq
\varnothing\mid S_0\bigl(\tau^{(0)}\bigr) \in J \bigr]
\cdot p(x_0) \\
& \geq&\mathbb{P} \bigl[\vphi\circ S_{0}\bigl[0,\tau^{(0)}\bigr]
\around^{(r)}
\xi\mid S_0\bigl(\tau^{(0)}\bigr) \in J \bigr] \cdot p(x_0).
\end{eqnarray*}
By weak convergence and Proposition~\ref{prop: mu-contin encoppasing a point},
and by Lemma~\ref{lem: BM from 0 several},
\[
\mathbb{P} \bigl[\vphi\circ S_{0}\bigl[0,\tau^{(0)}\bigr] \around
^{(r)} \xi\mid
S_0\bigl(\tau^{(0)}\bigr) \in J \bigr] \geq c_3,
\]
where $c_3 > 0$ is a constant.
Using (\ref{eqn: x0 at least 0}),
\[
p(0) \geq c_3 \cdot p(x_0) \geq c_4 \cdot p(w)
\]
for a constant $c_4 > 0$.
\end{pf*}

\subsection{\texorpdfstring{Proof of Proposition \protect\ref{prop: cond on hiting rhoj}}{Proof of Proposition 5.6}}
For $y \in V_\delta({D})$, define
\[
p(y) = \cases{
\mathbb{P} \bigl[S_y\bigl(\tau^{(y)}\bigr) = b \mid S_y\bigl(\tau
^{(y)}\bigr) \in J \bigr] ,
&\quad
if $\mathbb{P} \bigl[S_y\bigl(\tau^{(y)}\bigr) \in J \bigr] > 0$,
\cr
0, &\quad otherwise.}
\]
Since $p(\cdot)$ is harmonic, for every $y \in V_\delta(D)$,
there exists a path $\gamma_y$ from $y$ to $\p D$
such that $p(u) \geq p(y)$ for every $u \in\gamma_y$.
Let $w,y \in V_\delta(\rho_j)$.
Since $p(\cdot)$ is nonnegative, harmonic and bounded,
\[
p(w) \geq\mathbb{P} \bigl[S_w\bigl[0,\tau^{(w)}\bigr] \cap\gamma
_y \neq
\varnothing\mid S_w\bigl(\tau^{(w)}\bigr) \in J \bigr]
\cdot p(y).
\]
Let $\sigma^{(w)}$ be the first time $S_w$ exits $\vphi^{-1}(\rho
(\xi_j,\eta^2 R_j))$.
As in the proof of Lem\-ma~\ref{lem: starting near boundary},
if $\delta_0$ is small enough,
for every $v \sim u \in V_\delta(D)$, we have $|\vphi(v) - \vphi(u)|
< \beta(\eta- \eta^2)$.
By the strong Markov property,
\begin{eqnarray*}
&&\mathbb{P} \bigl[S_w\bigl[0,\tau^{(w)}\bigr] \cap\gamma_y \neq
\varnothing\mid
S_w\bigl(\tau^{(w)}\bigr) \in J \bigr]\\
&&\qquad= \frac{\mathbb{P} [S_w[0,\tau^{(w)}] \cap\gamma_y \neq
\varnothing, S_w(\tau^{(w)}) \in J ]}
{\mathbb{P} [S_w(\tau^{(w)}) \in J ]}
\\
&&\qquad\geq\frac{\mathbb{P} [S_w[0,\sigma^{(w)}] \cap\gamma_y
\neq\varnothing]
\cdot\min_{z} \mathbb{P} [S_z(\tau^{(z)}) \in J ]}
{\mathbb{P} [S_w(\tau^{(w)}) \in J ]} ,
\end{eqnarray*}
where the minimum is over $z \in V_\delta(D)$ such that $\vphi(z) \in
\rho(\xi_j,\eta R_j)$.
Define $h(z)$ to be the probability that $S_z(\tau^{(z)}) \in J$.
Since $h(\cdot)$ is harmonic,
there exists a path $g_w$ from $w$ to $\p D$
such that $h(u) \geq h(w)$ for every $u \in g_w$.
Since $h(\cdot)$ is nonnegative, harmonic and bounded,
by the choice of $\eta$,
\[
h(z) \geq\mathbb{P} \bigl[S_z\bigl[0,\tau^{(z)}\bigr] \cap g_w
\neq\varnothing
\bigr] \cdot p(w) \geq\biggl( 1 - \frac{\alpha\eps^2}{16}
\biggr) \cdot p(w).
\]
Also by the choice of $\eta$,
$\mathbb{P} [S_w[0,\sigma^{(w)}] \cap\gamma_y \neq\varnothing
] \geq1 - \frac{\alpha\eps^2}{16}$.
Thus,
\[
p(w) \geq\biggl( 1 - \frac{\alpha\eps^2}{8} \biggr) \cdot p(y).
\]
The strong Markov property implies the proposition.

\section{Convergence of the loop-erasure}
\label{sec: conv of lerw}

In this section we show that the scaling limit of the loop-erasure
of the reversal of the natural random walk on $G$ is $\SLE_2$
(for a planar-irreducible graph $G$ such that
the scaling limit of the natural random walk on $G$ is planar Brownian motion).
Most of our proof follows the proof of Lawler, Schramm and Werner in
\cite{LSW}.

\subsection{The observable}

Let $D \in\D$, and let $\delta> 0$.
For $v \in V_\delta(D)$,
let $S_v(\cdot)$ be the natural random walk on $G_\delta$
started at $v$ and stopped on exiting $D$.
Denote by $\hat{S}_v(\cdot)$ the loop-erasure of the reversal of
$S_v(\cdot)$.
\begin{rem}
There is a technicality we need to address.
Let $\gamma'(0),\ldots,\break\gamma'(T) = v$ be the loop-erasure\vadjust{\goodbreak}
of the reversal of $S_v(\cdot)$.
The edge $e = [\gamma'(0),\gamma'(1)]$ is not contained in $D$.
Define $\gamma(0) \in\p D$ as the last point on $e$ not in $D$
(see the definition of Poisson kernel in Section~\ref{sec: def and not}),
and define $\gamma(i) = \gamma'(i)$ for $i = 1,\ldots,T$.
\end{rem}

Let $\gamma(\cdot)$ be the loop-erasure of the reversal of a
natural random walk started at $0$ and stopped on exiting $D$;
that is, $\gamma(\cdot)$ has the same distribution as $\hat
{S_0}(\cdot)$,
but is independent of $S_0(\cdot)$
[from the time $\gamma(\cdot)$ hits $0$ it stays there].

\begin{prop} \label{prop: M_n martingale}
Let $v \in V_{\delta}(D)$. For $n \in\N$, define the random
variable
\[
M_n = \frac{ \mathbb{P} [ \hat{S}_v[0,n] = \gamma[0,n] ] }{ \mathbb
{P} [ \hat{S}_0[0,n] = \gamma[0,n] ] } .
\]
Then, $M_n$ is a martingale with respect to the filtration
generated by $\gamma[0,n]$.
\end{prop}
\begin{pf}
By the definition of $\gamma(\cdot)$, for every $w \in
V_{\delta}(D)$,
\[
\mathbb{P} \bigl[\gamma(n+1) = w \mid\gamma[0,n] \bigr] = \mathbb
{P} \bigl[\hat
{S_0}(n+1) = w \mid\hat{S_0}[0,n] = \gamma[0,n] \bigr] .
\]
Thus,
\begin{eqnarray*}
&&\E[M_{n+1} \mid\gamma[0,n] ] \\
&&\qquad= \sum_w \mathbb{P} \bigl[\gamma(n+1)
= w \mid\gamma[0,n] \bigr] \\
&&\qquad\quad\hspace*{11.3pt}{} \times\frac{ \mathbb{P} [\hat
{S_v}[0,n] =
\gamma[0,n] , \hat{S_v}(n+1) = w ] }{
\mathbb{P} [\hat{S_0}[0,n] = \gamma[0,n], \hat{S_0}(n+1)
= w ]}
\\
&&\qquad= \sum_w \mathbb{P} \bigl[\hat{S_v}(n+1) = w \mid\hat
{S_v}[0,n] =
\gamma[0,n] \bigr]
\frac{ \mathbb{P} [\hat{S_v}[0,n] = \gamma[0,n] ] }{
\mathbb{P} [\hat{S_0}[0,n] = \gamma[0,n] ] } \\
&&\qquad= M_n .
\end{eqnarray*}
\upqed
\end{pf}

Let $\Ee_n^{(v)}$ be the event that $S_v(\cdot)$
hits the set $\p D \cup\gamma[0,n]$ at $\gamma(n)$,
where we think of $S_v(\cdot)$ as a continuous curve
(linearly interpolated on the edges of $G_\delta$).
Denote
\[
H_n(v,\gamma(n)) = \mathbb{P} \bigl[ \Ee_n^{(v)} \bigr] .
\]

\begin{prop} \label{prop: H is a martingale}
For $v \in V_{\delta}(D)$,
\[
\frac{H_n(v, \gamma(n))}{ H_n(0, \gamma(n)) }
\]
is a martingale with respect to the filtration generated by
$\gamma[0,n]$.
\end{prop}
\begin{pf}
Define
\[
M_n = \frac{ \mathbb{P} [\hat{S_v}[0,n] = \gamma[0,n]
] }{ \mathbb{P} [\hat{S_0}[0,n] = \gamma[0,n] ] }\vadjust{\goodbreak}
\]
as in Proposition~\ref{prop: M_n martingale}. Since $M_n$ is a
martingale, it suffices to show that
\[
\frac{\mathbb{P} [ \Ee_n^{(v)} ]}{ \mathbb{P} [ \Ee_n^{(0)} ] } =
M_n .
\]

Let $z \in\{v,0\}$,
and let $S(\cdot)$ be the path $S_z[\T_z(\gamma[0,n]),\tau^{(z)}_D]$.
Since
$\{ \hat{S}_z[0,n] = \gamma[0,n] \} = \{ \hat{S}[0,n] = \gamma[0,n]
\}$,
by the strong Markov property,
\begin{eqnarray*}
\mathbb{P} \bigl[\hat{S_z}[0,n] = \gamma[0,n] , \Ee_n^{(z)}
\bigr] & = & \mathbb{P} \bigl[\hat{S}[0,n] = \gamma[0,n] ,
\Ee_n^{(z)} \bigr] \\
& = & \mathbb{P} \bigl[\hat{S}_{\gamma(n)}[0,n] = \gamma[0,n]
\bigr] \mathbb{P} \bigl[\Ee_n^{(z)} \bigr] ,
\end{eqnarray*}
which implies
%
%
\begin{equation} \label{eqn: Pr[Ee] Pr[s = g | Ee] first}
\mathbb{P} \bigl[\hat{S_z}[0,n] = \gamma[0,n] \mid\Ee_n^{(z)}
\bigr] =
\mathbb{P}
\bigl[\hat{S}_{\gamma(n)}[0,n] = \gamma[0,n] \bigr] .
\end{equation}

In addition, since $\{ \hat{S_z}[0,n] = \gamma[0,n] \} \subseteq\Ee
_n^{(z)}$,
%
%
\begin{equation} \label{eqn: Pr[Ee] Pr[s = g | Ee] second}
\mathbb{P} \bigl[\hat{S_z}[0,n] = \gamma[0,n] \bigr] =
\mathbb{P} \bigl[\Ee_n^{(z)} \bigr] \mathbb{P} \bigl[\hat
{S_z}[0,n] = \gamma
[0,n] \mid\Ee_n^{(z)} \bigr] .
\end{equation}

Combining (\ref{eqn: Pr[Ee] Pr[s = g | Ee] first}) and (\ref{eqn:
Pr[Ee] Pr[s = g | Ee] second}),
\begin{eqnarray*}
\frac{ \mathbb{P} [\Ee_n^{(v)} ] }{ \mathbb{P}
[\Ee_n^{(0)} ] }
&=& \frac{ \mathbb{P} [\Ee_n^{(v)} ] }{ \mathbb{P}
[\Ee_n^{(0)} ] } \cdot\frac{
\mathbb{P} [\hat{S}_{\gamma(n)} [0,n] = \gamma[0,n] ] }{
\mathbb{P} [\hat{S}_{\gamma(n)} [0,n] = \gamma[0,n] ]
} \\
&=& \frac{ \mathbb{P}
[\hat{S_v}[0,n] = \gamma[0,n] ] }{ \mathbb{P}
[\hat{S_0}[0,n] = \gamma[0,n] ] } = M_n .
\end{eqnarray*}
\upqed
\end{pf}

\subsection{The driving process}

Here are some known facts about the Schramm--Loewner evolution
(for more details, see~\cite{LSW}).
Let $D \in\D$, and let $\delta> 0$.
Let $\gamma(\cdot)$ be the loop-erasure of the reversal of a natural
random walk
started at $0$ and stopped on exiting $D$ (independent of $S_0$).
For $s \geq0$, define $\gamma[0,s]$ as the continuous curve
that is the linear interpolation of $\gamma(\cdot)$ on the edges of
$G_\delta$.
For $s \geq0$ such that $0 \notin\gamma[0,s]$,
define $\vphi_s\dvtx D \setminus\gamma[0,s] \to\U$ to be the unique
conformal map
satisfying $\vphi_s(0) = 0$ and $\vphi_s'(0) > 0$.
Let $t_s = \log\vphi_s'(0) - \log\vphi_D'(0)$, the \textit
{capacity} of
$\gamma[0,s]$ from $0$ in $D$. Let
\[
U_s = \lim_{z \to\gamma(s)} \vphi_s(z),
\]
where $z$ tends to $\gamma(s)$ from within $D \setminus\gamma[0,s]$.
Let $W\dvtx[0,\infty) \to\p\U$ be the unique continuous function such
that solving the
radial Loewner equation with driving function $W(\cdot)$ gives the
curve $\vphi_D \circ\gamma$.
Loewner's theory gives us the relation \mbox{$U_s = W(t_s)$}.
Let $\theta(\cdot)$ be the function such that $W(t) = W(0) e^{i
\theta(t)}$.
Let $\Delta_s = \theta(t_s)$, so we get that $U_s = U_0 e^{i
\Delta_s}$.
Since $t_s$ is a strictly increasing function of $s$,
we can define $\xi(r)$ to be the unique $s$ such that $t_{s} = r$
[by this definition, $\xi(t_r) = r$].
By the Loewner differential equation, for every $z \in D \setminus
\gamma[0,\xi(r)]$,
%
%
\begin{equation} \label{eqn: Loewner's diff eqn}
\p_r g_r(z) = g_r(z) \frac{U_{\xi(r)} + g_r(z)}{U_{\xi(r)} - g_r(z)
} ,
\end{equation}
where $g_r(z) = \vphi_{\xi(r)}(z)$.
\begin{prop} \label{prop: key estimate}
There exists $c > 0$ such that for all $\eps> 0$, there exists
$\delta_0 > 0$ such that for all $0 < \delta< \delta_0$
the following holds:

Let $D \in\D$.
Let $m = \min\{ 1 \leq j \in\N\dvtx t_j \geq\eps^2 \mbox{ or }
\vert\Delta_j\vert\geq\eps\} $. Then, a.s.,
\[
\vert\E[\Delta_m \mid\gamma(0) ] \vert\leq c \eps^3
\]
and
\[
\bigl\vert\E[\Delta_m^2 - 2 t_m \mid\gamma(0) ] \bigr\vert\leq
c \eps
^3 .
\]
\end{prop}
\begin{pf}
Fix $v \in V_{\delta}(D)$ such that $|\vphi_D(v)| \leq1/12$.
Let $Z = \vphi_0(v)$ and $U = U_0$.
We follow the proof of Proposition 3.4 in~\cite{LSW}, using our Lemma
\ref{lem: H(,) are close to lambda} (used with inner radius $c_1/8$)
to replace
Lemma 2.2 in~\cite{LSW}. This culminates to show that a.s.
%
%
\begin{eqnarray} \label{eqn: lambda_j - lambda_0}
&&\Ree\biggl(\frac{ZU(U+Z)}{(U-Z)^3} \biggr) \E[2 t_m - \Delta_m^2 \mid
\gamma(0) ]\nonumber\\[-8pt]\\[-8pt]
&&\qquad{}
+ \Imm\biggl(\frac{2ZU}{ (U-Z)^2} \biggr) \E[\Delta_m \mid\gamma(0) ]
= O(\eps^3) .\nonumber
\end{eqnarray}

Let $\eta= 1/20$.
Let $f(z) = \Ree(\frac{zU(U+z)}{(U-z)^3} )$ and $g(z) = \Imm(\frac
{2zU}{ (U-z)^2} )$.
We have
$f(\eta U) > 1/100$, $g(\eta U) = 0$, and $g(i \eta U) > 1/100$.
There exists $\eps' > 0$ such that for every
$z,w \in\frac{1}{12} \U$, if $|z-w| \leq\eps'$, then $|f(z)-f(w)|
\leq\eps^3$ and $|g(z)-g(w)| \leq\eps^3$.

Let $\D_{1,\eps'/2}$ be the finite family of domains given by
Proposition~\ref{prop: compact approx}.
By weak convergence,
there exists $\delta_0 > 0$ such that for all $0 < \delta< \delta_0$
and any $\tilde{D} \in\D_{1,\eps'/2}$,
there exist $v_1,v_2 \in V_{\delta}(\tilde{D})$ such that
$|\vphi_{\tilde{D}}(v_1) - \eta U| < \eps'/2$ and $|\vphi_{\tilde
{D}}(v_2) - i \eta U| < \eps'/2$.

Let $D \in\D$, and let $\tilde{D} \in\D_{1,\eps'/2}$ be the
$(1,\eps'/2)$-approximation of $D$.
Then, $\tilde{D} \subseteq D$ and $|\vphi_D(v_1) - \vphi_{\tilde
{D}}(v_1)| \leq\eps'/2$,
which implies that $f(\vphi_D(v_1)) = f(\eta U) + O(\eps^3)$ and
$g(\vphi_D(v_1)) = O(\eps^3)$.
Similarly, $g(\vphi_D(v_2)) = g(i \eta U) + O(\eps^3)$.
Applying (\ref{eqn: lambda_j - lambda_0}) to the vertices $v_1$ and
$v_2$, we have a.s.
\[
\bigl\vert\E[2 t_m - \Delta_m^2 \mid\gamma(0) ] \bigr\vert= O(\eps^3)
\quad\mbox{and}\quad \vert\E[\Delta_m \mid\gamma(0) ] \vert=
O(\eps^3) .
\]
\upqed\end{pf}

The following theorem shows that $\theta(\cdot)$
converges to one-dimensional Brownian motion.
\begin{theorem} \label{thm: driving process convergence}
For all $D \in\D$, and all $\alpha,T > 0$, there exists $\delta_0>
0$ such that for all $0< \delta< \delta_0$
the following holds:

Let $u \in[0,2\pi]$ be a uniformly distributed point, and let
$B_1(\cdot)$ be one-dimensional Brownian
motion started at $u$. Then, there is a coupling of $\gamma(\cdot)$
and $B_1(\cdot)$ such that
\[
\mathbb{P} \Bigl[\sup_{0 \leq t \leq T} | \theta(t) - B_1(2t) | >
\alpha\Bigr] < \alpha.
\]
\end{theorem}
\begin{pf}
The proof follows the proof of Theorem 3.7 in~\cite{LSW}, using our
Proposition~\ref{prop: key estimate} to replace Proposition 3.4 in
\cite{LSW}.
\end{pf}

\subsection{Weak convergence}

In this section we show that the scaling limit of the
loop-erasure of the reversal of the natural random walk on $G$ is $\SLE_2$.
It would seem natural to follow the proofs in Section 3.4 of~\cite{LSW}.
However, as stated in the \hyperref[sec1]{Introduction} there is a
difficulty with this approach.
The proof of tightness in~\cite{LSW} uses a ``natural'' family of
compact sets.
In our setting, it is not necessarily true that $\gamma$ belongs to
one of these
compact sets with high probability (and so the argument of~\cite{LSW} fails).
To overcome this difficulty, we define a ``weaker'' notion of tightness,
which we are able to use to conclude the proof.

\subsubsection{A sufficient condition for tightness}
\label{sec: tight}

For a metric space $\X$, and a set $A \subseteq\X$,
define $A^\eps= \bigcup_{a \in A} \rho(a,\eps)$, where
$\rho(a,\eps)$ is the ball of radius $\eps$ centered at $a$.
The following are Theorems 11.3.1, 11.3.3 and 11.5.4 in~\cite{Dudley}.

\begin{theorem}
Let $\X$ be a metric space. For any two laws $\mu,\nu$ on $\X$, let
\[
d(\mu,\nu) = \inf\{\eps> 0 \dvtx\mu(A) \leq\nu(A^\eps) +
\eps\mbox{ for all Borel sets } A \subset\X\} .
\]
Then, $d(\cdot,\cdot)$ is a metric on the space of laws on $\X$
[$d(\cdot,\cdot)$ is called the \textit{Prohorov metric}].
\end{theorem}
\begin{theorem} \label{thm: Prohorov metric}
Let $\X$ be a separable metric space. Let $\{\mu_n\}$
and $\mu$ be laws on $\X$.
Then, $\{\mu_n\}$ converges weakly to $\mu$ if and only
if $d(\mu_n,\mu) \to0$, where $d(\cdot,\cdot)$ is the
Prohorov metric.
\end{theorem}

Let $\{\mu_\delta\}$ be a family of laws on a metric
space $\X$.
We say that $\{\mu_{\delta}\}$ is \textit{tight} if for
every $\eps> 0$, there exists a compact set
$K_\eps\subset\X$ such that for all $\delta$, $\mu_{\delta
}(K_\eps) \geq1-\eps$.
\begin{theorem} \label{thm: Prohorov}
Let $\X$ be a complete separable metric space. Let $\{\mu
_\delta\}$ be a family of
laws on $\X$. Then, $\{\mu_\delta\}$ is tight if and
only if every sequence $\{\mu_{\delta_n}\}_{n \in\N}$
has a weakly-converging subsequence.
\end{theorem}

We use these theorems to prove an equivalent condition for tightness of measures
on a separable metric space.
\begin{lem} \label{lem: suff. cond. for tighness}
Let $\X$ be a complete separable metric space.
Let $\{\mu_m\}_{m \in\N}$ be a sequence of laws on $\X
$ with the following property:
for any $\eps> 0$, there exists a compact set $K_\eps\subset\X$
such that for any $\alpha> 0$, there exists $M > 0$ such that for all
$m \geq M$,
\[
\mu_m (K_{\eps}^{\alpha} ) \geq1-\eps.
\]
Then, the sequence $\{\mu_m\}$ is tight.
\end{lem}
\begin{pf}
Let $\{K_n\}$ be a sequence of compact sets such that for
all $\alpha> 0$,
there exists $M > 0$ such that for all $m \geq M$,
$\mu_m (K_{n}^{\alpha} ) \geq1- n^{-1}$.

Define
\[
M(\alpha,n) = \min\{j \in\N\dvtx\forall m \geq j \mu_m
(K_{n}^{\alpha} ) \geq1- n^{-1} \}.
\]
For $k \in\N$, define $M_0(1/k,n) = \max\{M(1/k,n), k \}$,
and for $\frac{1}{k} \leq\alpha< \frac{1}{k-1}$, define $M_0(\alpha
,n) = M_0(1/k,n)$.
For fixed $n$, the function $M_0(\cdot,n)$ has the following properties:
(i) The function $M_0(\alpha,n)$ is right-continuous in $\alpha$.
(ii) The function $M_0(\alpha,n)$ is a monotone nonincreasing
function of $\alpha$.
(iii) $\lim_{\alpha\to0} M_0(\alpha,n) = \infty$.
(iv) For every $0 < \alpha< 1$, $M_0(\alpha,n) \geq M(\alpha,n)$.

For every $m$, define $\alpha_n(m) = \inf\{ 0 < \beta< 1 \dvtx
M_0(\beta,n) \leq m \}$.
For every $\eta> 0$, $\alpha_n(M_0(\eta,n)) \leq\eta$, which
implies that
\[
\lim_{m \to\infty} \alpha_n(m) = 0 .
\]
In addition, $M_0(\alpha_n(m),n) \leq m$, which implies that for all
$m > 0$,
%
%
\begin{equation} \label{eqn: mu(m) of K(n,alpha)}
\mu_m \bigl(K_n^{\alpha_n(m)} \bigr) \geq1 - n^{-1} .
\end{equation}

For $m$ and $n \geq2$, define
\[
\mu_{m,n}(A) = \frac{\mu_m(A \cap K_n^{\alpha_n(m)})}{\mu
_m(K_n^{\alpha_n(m)}) }
\]
for all Borel $A \subset\X$.
We show that for any fixed $n \geq2$, the sequence $\{\mu
_{m,n}\}_{m \in\N}$ is tight.
Let $X_{m,n}$ be a random variable with law $\mu_{m,n}$.
Since $X_{m,n} \in K_n^{\alpha_n(m)}$ a.s., we can define a random variable
$\hat{X}_{m,n} \in K_n$ such that a.s.
the distance between $X_{m,n}$ and $\hat{X}_{m,n}$ is at most $2
\alpha_n(m)$.
Let $\hat{\mu}_{m,n}$ be the law of $\hat{X}_{m,n}$.
The Prohorov distance
between $\mu_{m,n}$ and $\hat{\mu}_{m,n}$ is at most $2 \alpha
_n(m)$. Thus,
if a sequence $\{\hat{\mu}_{m_k,n}\}_{k \in\N}$
converges to some limit in the Prohorov metric,
then the sequence $\{\mu_{m_k,n}\}_{k \in\N}$ has a
converging subsequence as well.
Since $\{\hat{\mu}_{m,n}\}$ is compactly supported,
it is a tight family of measures.
By Theorem~\ref{thm: Prohorov}, $\{\mu_{m,n}\}$ is also tight.

Thus, for any $n \geq2$ and any $\eps>0$, there
exists a compact set $K_{n,\eps} \subset\X$ such that for all $m >0$,
$\mu_{m,n}(K_{n,\eps}) \geq1-\eps$.
Let $\eps> 0$, and let $n = \lceil2/\eps\rceil$.
For all $m > 0$, by (\ref{eqn: mu(m) of K(n,alpha)}),
$\mu_{m} (K_n^{\alpha_n(m)}) \geq1-\eps/2$.
Thus,
\begin{eqnarray*}
\mu_m(K_{n,\eps/2}) &\geq&
\mu_m\bigl(K_{n,\eps/2} \cap K_n^{\alpha_n(m)}\bigr) =
\mu_{m,n}(K_{n,\eps/2}) \cdot\mu_m\bigl(K_n^{\alpha_n(m)}\bigr)\\
&\geq&(1-\eps/2)^2 > 1- \eps,
\end{eqnarray*}
which implies that the sequence $\{\mu_m\}$ is tight.
\end{pf}

\subsubsection{Quasi-loops}

Here we give some probability estimates needed for proving tightness.
\begin{clm}
\label{clm: z exits in zbeta}
Let $z \in\U$.
For all $\beta> 0$, there exist $c > 0$ and $\delta_0 > 0$ such that
for all $0 < \delta< \delta_0$
and for all $x \in V_{\delta}(\U)$ such that $|x-z| \geq2 \beta$,
\[
\mathbb{P} \bigl[S_x\bigl[0,\tau_{3 \U}^{(x)}\bigr] \cap\rho(z,\beta) =
\varnothing\bigr] \geq c .
\]
\end{clm}
\begin{pf}
It suffices to prove that there exists a set of vertices $U \subseteq
V_\delta(\U)$
such that every path starting at $x$ and reaching $\p\rho(x,\beta)$
intersects $U$,
and such that
\[
\mathbb{P} \bigl[S_u\bigl[0,\tau_{3 \U}^{(u)}\bigr] \cap\rho(z,\beta) =
\varnothing\bigr] \geq c
\]
for every $u \in U$.

Denote
$\mathcal{A} = \{ \frac{\beta}{100} (n + m \cdot i) \in\U\dvtx n,m
\in
\Z\}$.
The set $\mathcal{A}$ is finite, and there exists
${\tilde{x}} \in\mathcal{A}$ such that $x \in\rho({\tilde
{x}},\beta/40)$.

Assume toward a contradiction that such a set $U$ does not exist.
By the planarity of $G$,
there exists a path $Y \subseteq V_\delta(\U)$ in $G$
starting inside $\rho(\tilde x,\beta/40)$ and reaching $\p\rho
(\tilde x,\beta/2)$
such that
\[
\mathbb{P} \bigl[S_y\bigl[0,\tau_{3 \U}^{(y)}\bigr] \cap\rho(z,\beta) =
\varnothing\bigr] < c
\]
for every $y \in Y$.
On one hand, by weak convergence and Proposition~\ref{prop: mu-contin
encoppasing a point},
and by Proposition~\ref{prop: encompassing a point}
(and the conformal invariance of Brownian motion),
\begin{eqnarray*}
&&\mathbb{P} \bigl[ S_0\bigl[0,\tau^{(0)}_{3 \U}\bigr] \cap Y \neq\varnothing,
S_0\bigl[0,\tau^{(0)}_{3 \U}\bigr] \cap\rho(z,\beta) = \varnothing\bigr] \\
&&\qquad \geq\mathbb{P} \bigl[ S_0\bigl[0,\tau^{(0)}_{3 \U}\bigr] \around^{(\beta
/2)} \tilde x,
S_0\bigl[0,\tau^{(0)}_{3 \U}\bigr] \cap\rho(z,\beta) = \varnothing\bigr] > c .
\end{eqnarray*}
On the other hand,
\begin{eqnarray*}
&&\mathbb{P} \bigl[ S_0\bigl[0,\tau^{(0)}_{3 \U}\bigr] \cap Y \neq\varnothing,
S_0\bigl[0,\tau^{(0)}_{3 \U}\bigr] \cap\rho(z,\beta) = \varnothing\bigr] \\
&&\qquad \leq\max_{y \in Y}
\mathbb{P} \bigl[ S_y\bigl[0,\tau^{(y)}_{3 \U}\bigr] \cap\rho(z,\beta) =
\varnothing\bigr] < c,
\end{eqnarray*}
which is a contradiction.
%
\end{pf}
\begin{clm}
\label{clm: S intersect gamma}
There exist universal constants $c_1,c_2>0$ such that
for every $\eps> 0$
there exists $0 < C \leq c_1 \eps^{-c_2}$
such that for every $\beta>0$,
there exists $\delta_0 > 0$ such that for all
$0 < \delta< \delta_0$ the following holds:

Let $y \in V_{\delta}(\U)$ and let $g\dvtx[0,\infty] \to\C$ be a
curve such that
$g(0) \in\rho(y,\beta/C)$ and $g(\infty) \notin\rho(y,\beta)$.
Let $\tau_\beta$ be the exit time of $S_y(\cdot)$ from $\rho
(y,\beta)$. Then,
\[
\mathbb{P} \bigl[S_y [0, \tau_\beta] \cap g = \varnothing\bigr] <
\eps.
\]
\end{clm}
\begin{pf}
Let $c>0$ be the universal constant from Corollary~\ref{cor: w encompasses}
with the domain $2 \U$.
Let $N>1$ be large enough so that $(1-c)^N < \eps$, and let $C = 8
\cdot500^N$.
Denote
$\mathcal{A} = \{ \frac{\beta}{100 C} (n + m \cdot i) \in2 \U
\dvtx n,m \in\Z\}$.
There exists ${\tilde{y}} \in\mathcal{A}$ such that $y \in\rho
({\tilde{y}}, \frac{\beta}{40 C})$.
For $j=0,1,\ldots, N$,
let $r_j = 2 \cdot500^j \beta/ C$,
let $T_j$ be the first time $S_y(\cdot)$ exits $\rho(\tilde{y},400 r_j)$
and let $\Ee_j$ be the complement of the event $\{
S_y[T_j,T_{j+1}] \around^{(400r_{j+1})} \tilde{y} \}$.

By Corollary~\ref{cor: w encompasses}, there exists $\delta_0 > 0$
(independent of $y$, since $|\mathcal{A}| < \infty$) such that for
all $0 < \delta< \delta_0$,
we have $\mathbb{P} [ \Ee_0 ] \leq1-c$ and
$\mathbb{P} [ \Ee_j | \Ee_0,\ldots,\Ee_{j-1} ] \leq1-c$ for
all $j = 1,\ldots,N-1$.
Since $g$ is a continuous curve from $\rho(y,\beta/C) \subset\rho
(\tilde{y},2 \beta/C)$
to the exterior of $\rho(y,\beta) \supset\rho(\tilde{y},\beta/2)$,
and since $r_N < \beta/2$, for all $0 < \delta< \delta_0$,
\[
\mathbb{P} \bigl[S_y[0,\tau_\beta] \cap g = \varnothing\bigr] \leq
\mathbb{P} [\Ee_0,\Ee_1,\ldots,\Ee_{N-1} ] \leq
(1-c)^N < \eps.
\]
\upqed\end{pf}

Let $\gamma= \gamma_\delta$ be the loop-erasure of the reversal of
the natural random walk
on $V_{\delta}(\U)$, started at $0$ and stopped on exiting $\U$
($\gamma$ is a simple curve from $\p\U$ to $0$).
For $\alpha,\beta>0$, we say that $\gamma$ has a \textit{quasi-loop},
denoted $\gamma\in\QL(\alpha,\beta)$, if there exist $0 \leq s<
t<\infty$ such that
$|\gamma(s)- \gamma(t)| \leq\alpha$ and $\operatorname{diam}
(\gamma
[s,t]) \geq\beta$.

\begin{prop} \label{prop: quasi loop}
For all $\eps>0$ and all $\beta>0$, there exists $\alpha>0$ such
that for all $\delta>0$,
\[
\mathbb{P} [\gamma\in\QL(\alpha,\beta) ] < \eps.
\]
\end{prop}
\begin{pf}
Fix $\eps, \beta> 0$.
For $z \in\U$ and $\alpha> 0$,
let $\QL(z,\alpha,\beta)$ be the set of all curves $g$
such that there exist $0 \leq s < t < \infty$ such that
$g(s),g(t) \in\rho(z, \beta)$, $|g(s) - g(t)| \leq\alpha$
and $g[s,t] \not\subseteq\rho(z,2\beta)$.
Let $\mathcal{A} = \{ \frac{\beta}{100} (n + m \cdot i) \in\U
\dvtx n,m \in\Z\}$.

\begin{clm}
\label{clm: important for z}
For any $z \in\mathcal{A}$ and for any $\eta> 0$,
there exist $\alpha_1 > 0$ and $\delta_1 > 0$ such that for all
$0 < \delta< \delta_1$ the following holds:\vspace*{1pt}

Let $g$ be the loop-erasure of $S_0[0,\tau_{\U}^{(0)}]$
($g$ is \textup{not} the loop-erasure of the reversal). Then,
\[
\mathbb{P} [g \in\QL(z,\alpha_1,\beta) ] \leq\eta.
\]
\end{clm}
\begin{pf}
Fix $z \in\mathcal{A}$ and $\eta> 0$.
Let $s_1 \geq0$ be the first time $S_0(\cdot)$ hits $\rho(z,\beta)$,
and let $t_1 \geq s_1$ be the first time after $s_1$ that $S_0(\cdot)$
is not in $\rho(z,2\beta)$.
For $j \geq2$, let $s_j \geq t_{j-1}$ be the first time after
$t_{j-1}$ that $S_0(\cdot)$ hits $\rho(z,\beta)$,
and let $t_j \geq s_j$ be the\vadjust{\goodbreak} first time after $s_{j}$ that $S_0(\cdot
)$ is not in $\rho(z,2\beta)$.
Define $g_j$ as the loop-erasure of $S_0[0,t_j]$,
and let $Y_j$ be the event that $g_j \in\QL(z,\alpha_1,\beta)$.
Let $\tau= \tau^{(0)}_{3 \U}$,
and let $\mathcal{T}_j$ be the event that $t_j \leq\tau$.

Let $x$ be the first point on $g_j$ that is in $S_0[t_j,t_{j+1}]$.
Then, $g_{j+1}$ is $g_j$ up to the point $x$,
and then continues as the loop-erasure of $S_0[\sigma_x,t_{j+1}]$,
where $\sigma_x$ is the first time $S_0[t_j,t_{j+1}]$ hits $x$.

Denote $\mathcal{IO}_j = \{t_j \leq\tau< t_{j+1} \}$.
The event $\{s_{j} < \tau\}$ implies the event $\{
t_{j} < \tau\}$.
Thus, $\mathcal{IO}_j \cap\{g \in\QL(z,\alpha_1,\beta)
\} \subseteq Y_j ,$
which implies that for every $m \geq1$,
%
%
\begin{eqnarray}
\label{eqn: QL subset}
\{g \in\QL(z,\alpha_1,\beta) \} &\subseteq&
\mathcal{T}_{m+1} \cup\bigcup_{j = 1}^m \bigl(\{g \in\QL(z,\alpha
_1,\beta) \} \cap\mathcal{IO}_j\bigr) \nonumber\\[-8pt]\\[-8pt]
&\subseteq&
\mathcal{T}_{m+1} \cup\bigcup_{j = 1}^m Y_j.\nonumber
\end{eqnarray}

By Claim~\ref{clm: z exits in zbeta},
there exist $c > 0$ and $\delta_2 > 0$ such that for all $0 < \delta<
\delta_2$
and for all $x \in V_{\delta}(\U)$ such that $|x-z| \geq2 \beta$,
we have
$\mathbb{P} [S_x[0,\tau_{3 \U}^{(x)}] \cap\rho(z,\beta) =
\varnothing] \geq c$,
which implies that
%
%
\begin{equation}
\label{eqn: prob of tm}
\mathbb{P} [ \mathcal{T}_{m+1} ] \leq(1-c)^m < \eps/2
\end{equation}
for large enough $m$.

Fix $1 \leq j \leq m$. Let $h_{j+1}$ be the loop-erasure of $S_0[0,s_{j+1}]$.
Let $Q_j$ be the set of connected components of $h_{j+1} \cap\rho(z,2
\beta)$
that intersect $\rho(z,\beta)$ and are not connected to $S_0(s_{j+1})$.
By the definition of $s_{j+1}$, the size of $Q_j$ is at most $j$.

Assume that the event $Y_{j}$ does not occur.
If for every $K \in Q_j$, the distance between $S_0[s_{j+1},t_{j+1}]$
and $K \cap\rho(z,\beta)$ is more than $\alpha_1$,
then the event $Y_{j+1}$ does not occur.
Otherwise, let $K$ be the first component in $Q_j$ (according to the
order defined by time)
such that the distance between $S_0[s_{j+1},t_{j+1}]$ and $K \cap\rho
(z,\beta)$ is at most $\alpha_1$.
If $S_0[s_{j+1},t_{j+1}]$ intersects $K$,
then the event $Y_{j+1}$ does not occur.
Thus, the event $Y_{j+1} \setminus Y_j$ implies that there exists $K
\in Q_j$ such that
the distance between $S_0[s_{j+1},t_{j+1}]$ and $K \cap\rho(z,\beta
)$ is at most $\alpha_1$,
and $S_0[s_{j+1},t_{j+1}]$ does not intersect $K$.
By Claim~\ref{clm: S intersect gamma}, if $\alpha_1$ is small enough,
there exists $\delta_3 > 0$ such that for all $0 < \delta< \delta_3$,
since a.s. $|Q_j| \leq m$,
\[
\mathbb{P} [ Y_{j+1} \setminus Y_j ] < \frac{\eps}{2m} .
\]
Using (\ref{eqn: QL subset}) and (\ref{eqn: prob of tm}),
there exist $\alpha_1 > 0$ and $\delta_1 > 0$ such that for all \mbox{$0 <
\delta< \delta_1$},
\[
\mathbb{P} [ g \in\QL(z,\alpha_1,\beta) ] < \eps.
\]
\upqed\end{pf}

For every $z \in\U$, there exists
${\tilde{z}} \in\mathcal{A}$ such that $z \in\rho({\tilde
{z}},\beta/40)$.
Thus, for $\alpha< \beta/100$,
%
%
\begin{equation} \label{eqn: QL subset QL(z)}
\QL(\alpha,8 \beta) \subset\bigcup_{z \in\mathcal{A} } \QL(z,
\alpha, \beta) .
\end{equation}
Since the size of $\mathcal{A}$ does not depend on $\alpha$,
by Claim~\ref{clm: important for z},
there exist $\alpha_1 > 0$ and $\delta_1 > 0$ such that for all
$0 < \delta< \delta_1$,
and for every $z \in\mathcal{A}$,
\[
\mathbb{P} [g \in\QL(z,\alpha_1,\beta) ] < \frac
{\eps}{|\mathcal{A}|} ,
\]
which implies
\[
\mathbb{P} [g \in\QL(\alpha_1,8\beta) ] < \eps,
\]
where $g$ is the loop-erasure of $S_0[0,\tau_{\U}^{(0)}]$.

Let $\alpha_2 > 0$ be small enough so that for all
$z \in\mathcal{A}$ and all $\delta\geq\delta_1$,
we have that $\rho(z,\alpha_2)$ contains at most one vertex from
$G_\delta$.
Set $\alpha= \min\{\alpha_1, \alpha_2 \}$.
This implies that for any $\delta\geq\delta_1$,
$\mathbb{P} [g \in\QL(\alpha,8\beta) ] = 0$.
Therefore, for any $\delta> 0$,
\[
\mathbb{P} [g \in\QL(\alpha,8\beta) ] < \eps.
\]
By Lemma 1.1 in~\cite{Werner},
$g$ and $\gamma$ have the same law,
which completes the proof.
\end{pf}
\begin{prop} \label{prop: function for QL}
For every $\eps> 0$, there exists a monotone nondecreasing function
$f:(0,\infty) \to(0,1]$ such that
for all $\delta>0$,
\[
\mathbb{P} \bigl[\exists0 \leq s < t < \infty\dvtx\dist(\gamma
[0,s],\gamma[t,\infty]) < f( \diam(\gamma[s,t]) ) \bigr] < \eps.
\]
\end{prop}
\begin{pf}
By Proposition~\ref{prop: quasi loop}, for all $n \geq1$, there
exists $\alpha_n > 0$ such that for all $\delta>0$,
%
%
\begin{equation} \label{eqn: gamma has alpha(n) QL}
\sum_{n=1}^{\infty} \mathbb{P} [\gamma\in\QL(\alpha_n ,
2^{1-n} ) ] < \eps.
\end{equation}
Let $f\dvtx(0,\infty) \to(0,1]$ be a monotone nondecreasing function
such that
%
%
\begin{equation} \label{eqn: first condition on f}
f(2^{2-n}) < \alpha_n \qquad\mbox{for all } n \geq1 .
\end{equation}
Let $\delta> 0$.
Assume that there exist $0 \leq s < t < \infty$ such that
\[
\dist(\gamma[0,s],\gamma[t,\infty]) <
f( \diam(\gamma[s,t]) ) .
\]
Then, there exist $0 \leq s' < t' < \infty$ such that
$|\gamma(s') - \gamma(t')| < f ( \diam(\gamma[s',t']) )$.
Since $\gamma\subset\U$, there exists $n \geq1$ such that
$2^{1-n} < \diam(\gamma[s',t']) \leq2^{2-n}$.
By (\ref{eqn: first condition on f}),
there exists $n \geq1$ such that $|\gamma(s') - \gamma(t')| <
f(2^{2-n}) < \alpha_n$ and
$\diam(\gamma[s',t']) > 2^{1-n}$, which implies that
$\gamma\in\QL(\alpha_n, 2^{1-n})$.
The proposition follows by (\ref{eqn: gamma has alpha(n) QL}).
\end{pf}
\begin{prop} \label{prop: boundary problems far from boundary}
For every $\eps> 0$, there exists a monotone nondecreasing function
$f\dvtx(0,\infty) \to(0,1]$ such that
for every $\eta> 0$, there exists $\delta_0 > 0$ such that for every
$0 < \delta< \delta_0$,
\[
\mathbb{P} \bigl[\exists t \geq0 \dvtx\eta< 1- |\gamma(t)| < f(
\diam(\gamma[0,t]) ) \bigr] < \eps.
\]
\end{prop}
\begin{pf}
By Claim~\ref{clm: S intersect gamma} and the strong Markov property,
there exist universal constants $c_1,c_2>0$ such that
for every $m \geq1$,
there exists $0 < C_m \leq c_1 \eps^{-c_2} 2^{c_2 m}$ and $\delta_m >
0$ such that
for every $0 < \delta< \delta_m$,
%
%
\begin{equation} \label{eqn: small diameter}
\mathbb{P} \bigl[\diam\bigl( S_0\bigl[T(2^{1-m^2}),\tau_{\U}^{(0)}\bigr] \bigr) > C_m
2^{1-m^2} \bigr] < \eps2^{-m} ,
\end{equation}
where
\[
T(\xi) = \inf\{t \geq0 \dvtx1-|S_0(t)| \leq\xi\} .
\]
Since $C_m \cdot2^{-m^2}$ tends to $0$ as $m$ tends to infinity,
we can define a monotone nondecreasing function $f\dvtx(0,\infty) \to(0,1]$
such that $f(C_m 2^{1-m^2}) < 2^{1-(m+1)^2}$ for all $m \geq1$.

Denote by $\mathcal{Y}$
the event that there exists $t \geq0$ such that
$\eta< 1- |\gamma(t)| < f( \diam(\gamma[0,t]) )$.
Let $M$ be large enough so that $2^{1-M^2} < \eta$.
The event $\Y$ implies that
there exists $1 \leq m < M$ such that
\[
2^{1-(m+1)^2} < 1 - |\gamma(t)| \leq2^{1-m^2} ,
\]
which implies
\[
2^{1-(m+1)^2} < 1 - |\gamma(t)| < f(\diam(\gamma[0,t])) \leq
f \bigl(\diam\bigl(S_0\bigl[T(2^{1-m^2}),\tau_{\U}^{(0)}\bigr]\bigr) \bigr) .
\]
By the definition of $f$, this implies that $\diam(
S_0[T(2^{1-m^2}),\tau_{\U}^{(0)}] ) > C_m 2^{1-m^2}$.
Using (\ref{eqn: small diameter}),
for all $0 < \delta< \delta_0 = \min\{\delta_1,\ldots,\delta
_M\}$,
\[
\mathbb{P}[\mathcal{Y}] \leq\sum_{m = 1}^M \mathbb{P} \bigl[ \diam\bigl(
S_0\bigl[T(2^{1-m^2}),\tau_{\U}^{(0)}\bigr] \bigr) > C_m 2^{1-m^2} \bigr] < \eps.
\]
\upqed\end{pf}

\subsubsection{Tightness}

In this section we show that the laws of $\{\gamma_\delta
\}$ are tight.
Recall~$\mathcal{C}$, the space of all continuous curves with the
metric $\varrho$.
Let
\[
\X_0 = \{ g \in\mathcal{C} \dvtx g(0) \in\p\U, g(\infty) = 0 ,
g(0,\infty] \subset\U, g \mbox{ is a simple curve} \} .
\]
For a monotone nondecreasing function $f\dvtx(0,\infty) \to(0,1]$,
define $\X_f$ to be the set of $g \in\X_0$ such that for all $0 \leq
s < t < \infty$,
\[
\operatorname{dist}(g[0,s] \cup\p\U, g[t,\infty]) \geq f(
\operatorname{diam}(g[s,t]) ) .
\]
The following is Lemma 3.10 from~\cite{LSW}.
\begin{lem} \label{lem: compactness from LSW}
Let $f\dvtx(0,\infty) \to(0,1]$ be a monotone nondecreasing
function. Then,
$\X_f$ is compact in the topology of convergence with respect to the
metric~$\varrho$.
\end{lem}

For $\alpha>0$, define
\[
\X_f^{\alpha} = \{g \in\X_0 \dvtx\exists g' \in\X_f
\mbox{ such that } \varrho(g,g') < \alpha\} .
\]
\begin{lem} \label{lem: weak tightness of image of gamma}
For every $\eps>0$, there exists a monotone nondecreasing function
$f\dvtx(0,\infty) \to(0,1]$ such that
for any $\alpha> 0$, there exists $\delta_0>0$ such that for all $0 <
\delta<\delta_0$,
\[
\mathbb{P} [\gamma\notin\X_f^{\alpha} ] < \eps.
\]
\end{lem}
\begin{pf}
Let $g \in\X_0$ and let $\eta,\beta>0$.
Choose a parameterization for $g$ and let
$t_\eta(g) = \sup\{t \geq0 \dvtx1- |g(t)| \leq\eta\}
$. Define
$g^{\eta}$ to be the curve $g[t_\eta,\infty]$
(the curve $g^{\eta}$ does not depend on the choice of parameterization).
We say that a curve $h \in\X_0$ is $(\eta,\beta)$-\textit{adapted}
to $g$, if
$h^{\eta} = g^{\eta}$, and $\diam(h[0,t_\eta(h)] ) < \beta$.
Let $\A(g,\eta,\beta)$ be the set of all curves that are $(\eta
,\beta)$-adapted to $g$.
Note that $g$ is not necessarily in $\A(g,\eta,\beta)$, and that for
any two curves $h,\tilde{h} \in
\A(g,\eta,\beta)$,
%
%
\begin{equation} \label{eqn: small distance of Y}
\varrho(h,\tilde{h} ) \leq2 \beta.
\end{equation}

Define the curve $\tilde{\gamma}$ as follows.
Let $x \in\p(1-\eta)\U$ be the starting point of $\gamma^\eta$,
and let $y = \frac{x}{1-\eta} \in\p\U$.
Let $\tilde{\gamma}$ be the curve $[y,x] \cup\gamma^\eta$.

By Proposition~\ref{prop: boundary problems far from boundary},
there exists a monotone nondecreasing function $f_1\dvtx(0,\infty)
\to
(0,1]$ such that
for every $\eta> 0$, there exists $\delta_1 > 0$ such that for every
$0 < \delta< \delta_1$,
%
%
\begin{equation} \label{eqn: boundary problems}
\mathbb{P} \bigl[\exists t \geq0 \dvtx\eta< 1- |\gamma(t)| <
f_1( \diam(\gamma[0,t]) ) \bigr] < \eps/4 .
\end{equation}

By Proposition~\ref{prop: function for QL},
there exists a monotone nondecreasing function $f_2\dvtx(0,\infty)
\to
(0,1]$ such that
for all $\delta>0$,
%
%
\begin{equation} \label{eqn: QL}
\mathbb{P} \bigl[\exists0 \leq s < t < \infty\dvtx\dist(\gamma
[0,s],\gamma[t,\infty]) < f_2( \diam(\gamma[s,t]) ) \bigr] < \eps
/4 .\hspace*{-28pt}
\end{equation}

Define a monotone nondecreasing function $f\dvtx(0,\infty) \to(0,1]$ by
\[
f(\xi) = \min\{\xi/2,f_1(\xi/2),f_2(\xi/2)\}.
\]
Assume that there exists $t \geq0$ such that $1-|\tilde{\gamma}(t)|
< f(\diam(\tilde{\gamma}[0,t]))$.
Since $f(\xi) \leq\xi$,
there exists $t \geq0$ such that $\eta< 1- |\gamma(t)| < f( \diam
(\tilde{\gamma}[0,t]) )$,
and also $\diam( \tilde{\gamma}[0,t]) \leq\diam( \gamma[0,t] ) +
\eta$, which implies
\begin{eqnarray*}
\eta&<& f( \diam(\tilde{\gamma}[0,t]) )\leq\max\{f(2 \diam(
\gamma[0,t] )) , f(2 \eta) \}
\\ & \leq&\max\{f_1(\diam( \gamma[0,t] )) , \eta\}.
\end{eqnarray*}
Thus, there exists $t \geq0$ such that $\eta< 1- |\gamma(t)| < f_1(
\diam(\gamma[0,t]) )$.

Assume that there exist $0 \leq s < t < \infty$ such that
$|\tilde{\gamma}(t) - \tilde{\gamma}(s)| <\break f( \diam(\tilde{\gamma
}[s, t]) )$.
Let $t_\eta= t_\eta(\gamma)$.
Parameterize $\gamma$ and $\tilde{\gamma}$ so that $\gamma(t) =
\tilde{\gamma}(t)$
for every $t \geq t_\eta$.
Since $f(\xi) \leq\xi$, we have that $t > t_\eta$.
Assume that $s < t_\eta$.
Since $\diam(\tilde{\gamma}[s,t]) \leq\diam(\gamma[t_\eta,t]) +
|\tilde{\gamma}(t_\eta) - \tilde{\gamma}(s)|$,
\begin{eqnarray*}
|\tilde{\gamma}(t_\eta)- \tilde{\gamma}(s)| & \leq&
|\tilde{\gamma}(t) - \tilde{\gamma}(s)|
< f( \diam(\tilde{\gamma}[s,t]) ) \\
& \leq&\max\{f_2(\diam(\gamma[t_\eta,t])) , |\tilde{\gamma
}(t_\eta) - \tilde{\gamma}(s)| \} ,
\end{eqnarray*}
which implies
\[
|\gamma(t)- \gamma(t_\eta)| \leq|\tilde{\gamma}(t)- \tilde
{\gamma}(s)|
< f_2(\diam(\gamma[t_\eta,t])).
\]
If $s \geq t_\eta$, then
$|\gamma(t) - \gamma(s)| < f_2( \diam(\gamma[s,t]) )$.

Therefore, if $\tilde{\gamma} \notin\X_f$,
then either there exists $t \geq0$ such that $\eta< 1- |\gamma(t)| <
f_1( \diam(\gamma[0,t]) )$,
or there exist $0 \leq s < t < \infty$ such that
$|\gamma(t) - \gamma(s)| < f_2( \diam(\gamma[s,t]) )$.
By (\ref{eqn: boundary problems}) and (\ref{eqn: QL}),
for every $\eta> 0$,
there exists $\delta_1 > 0$ such that for every $0 < \delta< \delta_1$,
%
%
\begin{equation}
\label{eqn: gamma tilde in X}
\mathbb{P} [ \tilde{\gamma} \notin\X_f ] < \eps/2.
\end{equation}

By Claim~\ref{clm: S intersect gamma} and the strong Markov property,
for every $\alpha> 0$, there exist $\eta> 0$ and $\delta_2>0$ such
that for all $0< \delta< \delta_2$,
\[
\mathbb{P} \bigl[\diam\bigl(S_0\bigl[T(\eta),\tau_{\U}^{(0)}\bigr] \bigr) \geq\alpha
/2 \bigr] < \frac{\eps}{4} ,
\]
where $T(\eta) = \inf\{t \geq0 \dvtx1- |S_0(t)| \leq\eta
\}$.
If $1-|\gamma(t)| \leq\eta$, then
$\gamma[0,t] \subset S_0[T(\eta),\tau_{\U}^{(0)}]$.
Thus, for every $\alpha> 0$, there exist $0 < \eta< \alpha/2$ and
$\delta_2 > 0$ such that
for every $0 < \delta< \delta_2$,
%
%
\begin{equation} \label{eqn: gamma is in Y}
\mathbb{P} [ \varrho(\gamma,\tilde{\gamma}) \geq\alpha]
\leq\mathbb{P} [\gamma\notin\A(\gamma,\eta,\alpha/2)
] < \frac{\eps}{4} .
\end{equation}
Using (\ref{eqn: gamma tilde in X}),
for any $\alpha> 0$, there exist $\eta> 0$ and $\delta_0>0$ such
that for all \mbox{$0 < \delta<\delta_0$},
\[
\mathbb{P} [ \gamma\notin\X^\alpha_f ] \leq
\mathbb{P} [ \varrho(\gamma,\tilde{\gamma}) \geq\alpha] +
\mathbb{P} [ \tilde{\gamma} \notin\X_f ] < \eps.
\]
\upqed\end{pf}

Using Lemmas~\ref{lem: compactness from LSW},~\ref{lem: weak
tightness of image of gamma}
and~\ref{lem: suff. cond. for tighness},
we have the following corollary.
\begin{cor} \label{cor: tightness of image of gamma}
Let $\{\delta_n\}$ be a sequence converging to zero,
and let $\mu_{n}$ be the law of the curve $\gamma_{\delta_n}$.
Then, the sequence $\{\mu_n\}$ is tight.
\end{cor}

\subsubsection{Convergence}

Here we finally show that the scaling limit of the
loop-erasure of the reversal of the natural random walk on $G$ is $\SLE_2$.
We first show that any subsequential limit of $\{\gamma_\delta
\}$
is a.s. a simple curve.
\begin{lem} \label{lem: simple curve}
Let $\{\delta_n\}$ be a sequence converging to zero,
and let $\mu_{n}$ be the law of the curve $\gamma_{\delta_n}$.
If $\mu_n$ converges weakly to $\mu$,
then $\mu$ is supported on~$\X_0$.
\end{lem}
\begin{pf}
Let $d(\cdot,\cdot)$ be the Prohorov metric.
By Theorem~\ref{thm: Prohorov metric}, \mbox{$d(\mu_n,\mu) \to0$}.

As in the proof of Lemma~\ref{lem: weak tightness of image of gamma},
by (\ref{eqn: gamma tilde in X}) and (\ref{eqn: gamma is in Y}),
for every $\eps> 0$, there exists a monotone nondecreasing function
$f\dvtx(0,\infty) \to(0,1]$
such that for every $\alpha> 0$,
there exists $\delta_0 > 0$ such that for all $0< \delta< \delta_0$,
we can define a curve $\gamma_{\delta}^\alpha$ such that
%
%
\begin{equation} \label{eqn: gamma(alpha)}
\mathbb{P} [ \gamma_\delta^{\alpha} \notin\X_f ] < \eps
\quad\mbox{and}\quad
\mathbb{P} [ \varrho(\gamma_\delta, \gamma_\delta^\alpha) \geq
\alpha] < \alpha.\vadjust{\goodbreak}
\end{equation}
Let $\mu_n^{\alpha}$ be the law of $\gamma_{\delta_n}^{\alpha}$.
By (\ref{eqn: gamma(alpha)}),
for all $k \in\N$, there exists $f_k$ such that
for every $m \in\N$, there exists $N_{m,k} > m+k$ such that for all
$n \geq N_{m,k}$,
we have $d(\mu_n,\mu_n^{1/m}) < 1/m$ and $\mu_n^{1/m}(\X_{f_k}) > 1-1/k$.

Since $d(\mu_{N_{m,k}}^{1/m},\mu) \leq d(\mu_{N_{m,k}}^{1/m},\mu
_{N_{m,k}}) + d(\mu_{N_{m,k}},\mu)$,
by Theorem~\ref{thm: Prohorov metric},
for every fixed $k \in\N$,
the sequence $\{ \mu_{N_{m,k}}^{1/m} \}_{m \in\N}$
converges weakly\vspace*{1pt} to $\mu$.
Using Lemma~\ref{lem: compactness from LSW},
the Portmanteau theorem (see Chapter III in~\cite{Shiry}) tells us
that for every $k \in\N$,
\[
\mu(\X_{f_k}) \geq\mathop{\lim\sup}_{m \to\infty} \mu
_{N_{m,k}}^{1/m} (\X_{f_k}) > 1- 1/k .
\]
Thus, since $\X_{f_k} \subseteq\X_0$ for all $k \in\N$,
\[
\mu(\X_0) \geq\mu\biggl( \bigcup_k \X_{f_k} \biggr) = 1 .
\]
\upqed\end{pf}
\begin{pf*}{Proof of Theorem~\ref{thm: main thm}}
The proof follows by plugging
Theorem~\ref{thm: driving process convergence}, Corollary~\ref{cor:
tightness of image of gamma},
and Lemma~\ref{lem: simple curve}
into the proof of Theorem 3.9 in~\cite{LSW}.
\end{pf*}

\section*{Acknowledgments}
We thank Itai Benjamini for suggesting the problem,
and for very helpful discussions. We also thank Gady Kozma for helpful
discussions, and Ofer Zeitouni and Nathanael Berestycki for their help
with Lemma~\ref{lem: suff. cond. for tighness}.

%

%
\printaddresses


\begin{thebibliography}{17}

\bibitem{BerBisk}
%
\begin{barticle}[mr]
\bauthor{\bsnm{Berger},~\bfnm{Noam}\binits{N.}} \AND
\bauthor{\bsnm{Biskup},~\bfnm{Marek}\binits{M.}}
(\byear{2007}).
\btitle{Quenched invariance principle for simple random walk on percolation
clusters}.
\bjournal{Probab. Theory Related Fields}
\bvolume{137}
\bpages{83--120}.
\bid{doi={10.1007/s00440-006-0498-z}, mr={2278453}}
\end{barticle}
%
\endbibitem

\bibitem{Conway}
%
\begin{bbook}[mr]
\bauthor{\bsnm{Conway},~\bfnm{John~B.}\binits{J.~B.}}
(\byear{1995}).
\btitle{Functions of One Complex Variable. {II}}.
\bseries{Graduate Texts in Mathematics}
\bvolume{159}.
\bpublisher{Springer}, \baddress{New York}.
\bid{mr={1344449}}
\end{bbook}
%
\endbibitem

\bibitem{Dudley}
%
\begin{bbook}[mr]
\bauthor{\bsnm{Dudley},~\bfnm{Richard~M.}\binits{R.~M.}}
(\byear{1989}).
\btitle{Real Analysis and Probability}.
\bpublisher{Wadsworth \& Brooks/Cole Advanced Books \& Software},
\baddress{Pacific Grove, CA}.
\bid{mr={0982264}}
\end{bbook}
%
\endbibitem

\bibitem{Kozma}
%
\begin{barticle}[mr]
\bauthor{\bsnm{Kozma},~\bfnm{Gady}\binits{G.}}
(\byear{2007}).
\btitle{The scaling limit of loop-erased random walk in three dimensions}.
\bjournal{Acta Math.}
\bvolume{199}
\bpages{29--152}.
\bid{doi={10.1007/s11511-007-0018-8}, mr={2350070}}
\end{barticle}
%
\endbibitem

\bibitem{LERWdfn}
%
\begin{barticle}[mr]
\bauthor{\bsnm{Lawler},~\bfnm{Gregory~F.}\binits{G.~F.}}
(\byear{1980}).
\btitle{A self-avoiding random walk}.
\bjournal{Duke Math. J.}
\bvolume{47}
\bpages{655--693}.
\bid{mr={0587173}}
\end{barticle}
%
\endbibitem

\bibitem{LawRWRE}
%
\begin{barticle}[mr]
\bauthor{\bsnm{Lawler},~\bfnm{Gregory~F.}\binits{G.~F.}}
(\byear{1982}).
\btitle{Weak convergence of a random walk in a random environment}.
\bjournal{Comm. Math. Phys.}
\bvolume{87}
\bpages{81--87}.
\bid{mr={0680649}}
\bptnote{check year}
\end{barticle}
%
\endbibitem

\bibitem{LawlerSurvey}
%
\begin{bincollection}[mr]
\bauthor{\bsnm{Lawler},~\bfnm{Gregory~F.}\binits{G.~F.}}
(\byear{1999}).
\btitle{Loop-erased random walk}.
In \bbooktitle{Perplexing Problems in Probability}.
\bseries{Progress in Probability}
\bvolume{44}
\bpages{197--217}.
\bpublisher{Birkh\"auser}, \baddress{Boston, MA}.
\bid{mr={1703133}}
\end{bincollection}
%
\endbibitem

\bibitem{LawlerBook}
%
\begin{bbook}[mr]
\bauthor{\bsnm{Lawler},~\bfnm{Gregory~F.}\binits{G.~F.}}
(\byear{2005}).
\btitle{Conformally Invariant Processes in the Plane}.
\bseries{Mathematical Surveys and Monographs}
\bvolume{114}.
\bpublisher{Amer. Math. Soc.}, \baddress{Providence, RI}.
\bid{mr={2129588}}
\end{bbook}
%
\endbibitem

\bibitem{LSW}
%
\begin{barticle}[mr]
\bauthor{\bsnm{Lawler},~\bfnm{Gregory~F.}\binits{G.~F.}},
\bauthor{\bsnm{Schramm},~\bfnm{Oded}\binits{O.}} \AND
\bauthor{\bsnm{Werner},~\bfnm{Wendelin}\binits{W.}}
(\byear{2004}).
\btitle{Conformal invariance of planar loop-erased random walks and uniform
spanning trees}.
\bjournal{Ann. Probab.}
\bvolume{32}
\bpages{939--995}.
\bid{doi={10.1214/aop/1079021469}, mr={2044671}}
\end{barticle}
%
\endbibitem

\bibitem{PeresBM}
%
\begin{bbook}[auto:STB|2010-11-18|09:18:59]
\bauthor{\bsnm{M{\"o}rters},~\bfnm{P.}\binits{P.}} \AND
\bauthor{\bsnm{Peres},~\bfnm{Y.}\binits{Y.}}
(\byear{2010}).
\btitle{Brownian Motion}.
\bpublisher{Cambridge Univ. Press}, \baddress{Cambridge}.
\end{bbook}
%
\endbibitem

\bibitem{Pom}
%
\begin{bbook}[mr]
\bauthor{\bsnm{Pommerenke},~\bfnm{Ch.}\binits{C.}}
(\byear{1992}).
\btitle{Boundary Behaviour of Conformal Maps}.
\bseries{Grundlehren der Mathematischen Wissenschaften [Fundamental Principles
of Mathematical Sciences]}
\bvolume{299}.
\bpublisher{Springer}, \baddress{Berlin}.
\bid{mr={1217706}}
\end{bbook}
%
\endbibitem


\bibitem{RS01}
%
\begin{barticle}[mr]
\bauthor{\bsnm{Rohde},~\bfnm{Steffen}\binits{S.}} \AND
\bauthor{\bsnm{Schramm},~\bfnm{Oded}\binits{O.}}
(\byear{2005}).
\btitle{Basic properties of {SLE}}.
\bjournal{Ann. of Math. (2)}
\bvolume{161}
\bpages{883--924}.
\bid{doi={10.4007/annals.2005.161.883}, mr={2153402}}
\end{barticle}
%
\endbibitem

\bibitem{SchSLE}
%
\begin{barticle}[mr]
\bauthor{\bsnm{Schramm},~\bfnm{Oded}\binits{O.}}
(\byear{2000}).
\btitle{Scaling limits of loop-erased random walks and uniform
spanning trees}.
\bjournal{Israel J. Math.}
\bvolume{118}
\bpages{221--288}.
\bid{doi={10.1007/BF02803524}, mr={1776084}}
\end{barticle}
%
\endbibitem

\bibitem{Shiry}
%
\begin{bbook}[mr]
\bauthor{\bsnm{Shiryaev},~\bfnm{A.~N.}\binits{A.~N.}}
(\byear{1996}).
\btitle{Probability},
\bedition{2nd} ed.
\bseries{Graduate Texts in Mathematics}
\bvolume{95}.
\bpublisher{Springer}, \baddress{New York}.
\bid{mr={1368405}}
\end{bbook}
%
\endbibitem

\bibitem{Werner}
%
\begin{bincollection}[mr]
\bauthor{\bsnm{Werner},~\bfnm{Wendelin}\binits{W.}}
(\byear{2004}).
\btitle{Random planar curves and {S}chramm--{L}oewner evolutions}.
In \bbooktitle{Lectures on Probability Theory and Statistics}.
\bseries{Lecture Notes in Math.}
\bvolume{1840}
\bpages{107--195}.
\bpublisher{Springer}, \baddress{Berlin}.
\bid{mr={2079672}}
\bptnote{check year}
\end{bincollection}
%
\endbibitem

\end{thebibliography}
\end{document}